\documentclass[10pt,a4paper]{article}

\usepackage{amsmath, amsfonts, amssymb, amsthm}
\usepackage{mathrsfs}
\usepackage{graphicx}
\usepackage{geometry}
\usepackage{float}
\usepackage{fancyhdr}
\usepackage{siunitx}
\usepackage{listings}
\usepackage{caption}
\usepackage{booktabs}
\usepackage{multirow}
\usepackage{array}
\usepackage{adjustbox}
\usepackage{algorithm}
\usepackage{algpseudocode}
\usepackage[dvipsnames]{xcolor}
\usepackage[skip = \baselineskip, indent = 0pt]{parskip}
\usepackage{parskip}
\usepackage{lastpage}
\usepackage{tikz}
\usepackage{pgfplots}
\usepackage{subcaption}
\usepackage[OT2, T1]{fontenc}
\usepackage{fancyvrb}
\usepackage{blkarray}
\usepackage{mathtools}
\usepackage{xparse}
\usepackage{epigraph}
\usepackage{afterpage}
\usepackage{lettrine}
\usepackage{lmodern}
\usepackage{tabularx}
\usepackage{bookmark}
\usepackage[
  backend = biber,
  style = numeric,
  sorting = ynt
]{biblatex}
\usepackage{enumerate}
\usepackage{hyperref}
\usepackage[capitalize, noabbrev]{cleveref}
\usepackage{alphalph}
\usepackage{authblk}

\pgfplotsset{compat = newest}
\usepgfplotslibrary{statistics}
\usepgfplotslibrary{groupplots}

\DeclareSymbolFont{cyrletters}{OT2}{wncyr}{m}{n}
\DeclareMathSymbol{\Sha}{\mathalpha}{cyrletters}{"58}

\usetikzlibrary{shapes, arrows, positioning, calc}

\tikzstyle{block} = [draw, rectangle, minimum width = 3em, minimum height = 3em]
\tikzstyle{sum} = [draw, circle, node distance = 1.5cm]
\tikzstyle{mult} = [draw, regular polygon, regular polygon sides = 3, shape border rotate = -90, node distance = 1cm, inner sep = 1mm, outer sep = 0mm, text width = 1em]
\tikzstyle{multiso} = [draw, isosceles triangle, node distance = 1cm, inner sep = 1mm, outer sep = 0mm, text width = 1em]
\tikzstyle{input} = [coordinate]
\tikzstyle{output} = [coordinate]

\makeatletter
\pgfplotsset{
  auto title/.style={
    title = (\alphalph{\pgfplots@group@current@plot}),
    font = \scriptsize
  }
}
\makeatother

\makeatletter
\def\thickhrulefill{\leavevmode \leaders \hrule height 1ex \hfill \kern \z@}
\def\@makechapterhead#1{%
  \vspace*{10\p@}%
  {\parindent \z@ \centering \reset@font
        \thickhrulefill\quad
        {\scshape \@chapapp{} \thechapter}
        \quad \thickhrulefill
        \par\nobreak
        \vspace*{10\p@}%
        \interlinepenalty\@M
        \hrule
        \vspace*{10\p@}%
        \Huge \bfseries #1\par\nobreak
        \par
        \vspace*{10\p@}%
        \hrule
    \vskip 100\p@
  }}
\def\@makeschapterhead#1{%
  \vspace*{10\p@}%
  {\parindent \z@ \centering \reset@font
        \thickhrulefill
        \par\nobreak
        \vspace*{10\p@}%
        \interlinepenalty\@M
        \hrule
        \vspace*{10\p@}%
        \Huge \bfseries #1\par\nobreak
        \par
        \vspace*{10\p@}%
        \hrule
    \vskip 100\p@
  }}

  \makeatletter

  \newcommand\backmatter{%
    \if@openright
      \cleardoublepage
    \else
      \clearpage
    \fi
     }
  
\makeatother

\makeatletter
\newenvironment{breakablealgorithm}
  {
   \begin{center}
     \refstepcounter{algorithm}
     \hrule height.8pt depth0pt \kern2pt
     \renewcommand{\caption}[2][\relax]{
       {\raggedright\textbf{\fname@algorithm~\thealgorithm} ##2\par}%
       \ifx\relax##1\relax 
         \addcontentsline{loa}{algorithm}{\protect\numberline{\thealgorithm}##2}%
       \else 
         \addcontentsline{loa}{algorithm}{\protect\numberline{\thealgorithm}##1}%
       \fi
       \kern2pt\hrule\kern2pt
     }
  }{
     \kern2pt\hrule\relax
   \end{center}
  }
\makeatother

\hypersetup{
    pdfencoding = auto,
    psdextra,
    colorlinks = true,
    linkcolor = Magenta,
    filecolor = Green,
    urlcolor = Cyan,
    citecolor = Cyan,
    pdftitle = {GSLRA},
    pdfauthor = {Rishav Chakraborty}
}

\newcommand\numberthis{\addtocounter{equation}{1}\tag{\theequation}}

\theoremstyle{definition}
\newtheorem{definition}{Definition}[section]
\newtheorem{notation}[definition]{Notation}

\newtheorem{remark}[definition]{Remark}
\theoremstyle{plain}
\newtheorem{theorem}[definition]{Theorem}

\newtheorem{lemma}[definition]{Lemma}

\crefrangeformat{axiom}{axioms #3#1#4 to #5#2#6}
\crefalias{problem}{equation}
\crefname{problem}{problem}{problems}
\Crefname{problem}{Problem}{Problems}
\creflabelformat{problem}{#2\textup{(#1)}#3}
\crefname{theorem}{Theorem}{Theorems}

\algdef{SE}{Begin}{End}{\textbf{Begin}}{\textbf{End}}

\newcommand*{\N}{\mathbb{N}}
\newcommand*{\NZ}{\mathbb{N}_0}
\newcommand*{\Z}{\mathbb{Z}}

\newcommand*{\Q}{\mathbb{Q}}
\newcommand*{\R}{\mathbb{R}}

\renewcommand*{\C}{\mathbb{C}}

\newcommand*{\Order}{\mathcal{O}}

\newcommand*{\inv}{^{-1}}

\newcommand*{\trans}{^\mathsf{T}}

\renewcommand*{\L}{\mathcal{L}}
\renewcommand*{\P}{\mathscr{P}}
\renewcommand*{\Pr}{\mathbb{P}}
\newcommand*{\E}{\mathbb{E}}

\newcommand*{\pt}{\mathcal{P}}
\newcommand*{\defeq}{\coloneqq}

\newcommand*{\norm}[2][\DefaultNormSubscript]{\left\Vert {#2} \right\Vert_{#1}}

\newcommand{\innerproduct}[3][\DefaultInnerProductSubscript]{\left\langle {#2}, {#3} \right\rangle_{#1}}
\newcommand*{\Fr}{\mathrm{F}}

\newcommand*{\irng}[1]{\mathfrak{I}_{#1}}
\def\DefaultIndex{j}
\def\DefaultLower{1}
\def\DefaultUpper{\infty}
\def\AddSubscript{1}
\def\DefaultAddSubscript{\AddSubscript}
\NewDocumentCommand{\seq}{ O{\DefaultIndex} O{\DefaultAddSubscript} O{\DefaultLower} O{\DefaultUpper} m }{\left\{ \if#2\AddSubscript {#5}_{#1} \else {#5} \fi \right\}_{#1 = #3}^{#4}}
\def\OptSymbol{*}
\newcommand*{\opt}{^\OptSymbol}
\newcommand*{\optTrans}{^{\OptSymbol\mathsf{T}}}

\DeclareMathOperator*{\supp}{supp}

\DeclareMathOperator*{\Null}{Null}
\DeclareMathOperator*{\rank}{rank}
\DeclareMathOperator*{\nullity}{nullity}

\DeclareMathOperator{\vectorize}{vec}
\DeclareMathOperator{\Vect}{vec}
\DeclareMathOperator{\spec}{spec}

\begin{document}
\title{\LARGE \bf A Generalized Framework for Structured Low-Rank Approximation and Completion}
\author[1]{Rishav Chakraborty}
\author[1]{A. Sanand Amita Dilip}
\author[2]{Swanand R. Khare}
\affil[1]{{\small Department of Electrical Engineering, Indian Institute of Technology Kharagpur, Kharagpur, India.}}
\affil[2]{{\small Department of Mathematics, Indian Institute of Technology Kharagpur, Kharagpur, India.}}
\date{}
\maketitle
\begin{abstract}
    We formulate a generalized version of SLRA, consisting of general affine constraints, which subsumes the problems of structured low-rank matrix approximation and completion as special cases. We characterize its critical points and propose an algorithm to approximately solve the general problem. Finally, we demonstrate the performance of the proposed algorithm for stock-price forecasting, computing approximate gcd of polynomials with noisy coefficients, and time-series denoising.
\end{abstract}
\section{Introduction}
\label{sec:Introduction}
Low-rank approximation is a foundational tool in modern data analysis, arising naturally whenever high-dimensional data is generated by, or well-approximated by, a small number of latent factors \cite{haeffele-vidal-2020,ishteva-usevich-markovsky, gillard-usevich}. Given a data matrix, the goal is to find a matrix of low rank that best represents the observed data under some notion of fidelity. 
This simple idea underlies techniques as diverse as principal component analysis, system identification, signal denoising, collaborative filtering, and dimensionality reduction, and it has become a central object of study in numerical linear algebra, statistics, and machine learning alike.

In many applications, however, the matrix of interest is constrained to possess additional algebraic structure. 
Hankel and Toeplitz matrices arise in time-series analysis \cite{gillard-usevich,reddi}, system identification \cite{fazel-pong-sun-tseng,wahlberg}, and signal processing \cite{zhang-zhang,strang}, where they encode shift-invariance in the underlying dynamical process. Sylvester \cite{zhi-approximate-gcd} and B\'ezout matrices \cite{bitmead-kung-anderson-kailath} appear in computer algebra and control theory, where they encode polynomial relationships. Covariance matrices in statistical estimation must be symmetric and positive semi-definite. In each of these settings, an approximation that ignores the structural constraints---treating the problem as a generic unstructured low-rank approximation---can produce a solution that is numerically close to the data but structurally meaningless, and the loss of structure can propagate into downstream errors that are far more consequential than the fitting error itself. This has motivated a substantial body of work on structured low-rank approximation (SLRA): the problem of finding, within a prescribed linear (or affine) family of structured matrices, one whose rank does not exceed a given bound and which is as close as possible, in an appropriate norm, to a given data matrix. Unlike its unstructured counterpart, which admits a closed-form solution via the singular value decomposition (due to the Eckart-Young theorem), SLRA is in general NP-hard, non-convex, and lacks a closed-form solution, which has led to a rich landscape of algorithmic approaches ranging from alternating projections and variable projection methods to semi-definite relaxations \cite{cifuentes} and Riemannian optimization \cite{vandereycken}. The simplest of the family of alternating projections algorithms is the Cadzow iteration \cite{cadzow}, which repeatedly projects the matrix in the space of all matrices given by a structure and the space of all matrices with the given rank in an alternating fashion. It was shown in \cite{knirsch-petz-plonka} that this algorithm always converges for rank \(1\) matrices, and for rank \(r\) matrices with \(r\in\N\), the sequence of iterates always contains a convergent subsequence \cite{gillard-zhigljavsky-2013,zvonarev-golyandina}. There exist specialized algorithms for approximately solving SLRA for specific structures such as multi-start APBR \cite{gillard-zhigljavsky}, and variable projections \cite{markovsky-usevich} for Hankel structures. Furthermore, some general algorithms have had great success when applied to specific structures such as STLN applied for Sylvester and generalized Sylvester matrices \cite{zhi-approximate-gcd}.

A closely related but distinct problem is matrix completion: recovering an entire low-rank matrix from only a sparse subset of its entries. Popularized by applications such as recommender systems---most famously the Netflix Prize---matrix completion asks under what conditions, and via what algorithms, a low-rank matrix can be exactly or approximately recovered when the vast majority of its entries are unobserved. The seminal work of Candes \cite{candes-plan-yaniv,candes-recht}, Recht \cite{recht-fazel-parrilo}, Tao \cite{candes-tao}, and others established that, under incoherence conditions on the singular vectors and a sufficient number of randomly sampled entries, nuclear norm minimization recovers the underlying low-rank matrix with high probability, and this convex relaxation perspective has since been extended to noisy, structured, and non-uniformly sampled settings \cite{candes-plan,chen-chi-fan-ma-yan,negahban-wainwright-2011,negahban-wainwright-2012}. Viewed through the right lens, matrix completion can itself be seen as a special case of structured low-rank approximation, where the ``structure'' is the pattern of observed and missing entries, and the approximation criterion becomes a constraint on the observed positions. Conversely, many structured low-rank approximation problems---particularly those involving Hankel or Toeplitz structure---can be reformulated as completion problems on an associated structured matrix, since imposing exact structure is equivalent to requiring consistency across repeated entries. This interplay between the two problems suggests that algorithmic and theoretical advances in one domain can often be transferred to the other, and yet the literature on SLRA and on matrix completion has, until relatively recently, developed largely in parallel, with limited cross-pollination of techniques.

This paper aims to bridge these two threads.
We consider a general version of low-rank approximation of real matrices with affine equality constraints as shown in (\ref{eqn:TheGeneralOptimizationProblem}), and analyze the curvature of its constrained critical points. Two special forms of such a formulation are structured low-rank approximation (SLRA) and structured low-rank matrix completion (SLRC). Finally, we propose an algorithm to obtain a ``good'' local minimum of such an optimization problem. In particular, we demonstrate results on the SLRA and SLRC of Hankel matrices, and SLRA of Sylvester matrices. 


Our contributions are three-fold:
\begin{enumerate}
    \item We formulate a generalized version of structured low-rank matrix approximation problem as shown in (\ref{eqn:TheGeneralOptimizationProblem}) that subsumes some important classes of optimization problems, namely, SLRA and SLRC.
    \item We characterize (see \cref{thm:CriticalPoints}) the constrained critical points of the formulated optimization problem.
    \item We propose an algorithm to find a ``good'' local minimum of (\ref{eqn:TheGeneralOptimizationProblem}), and validate it on the problems of stock-price forecasting (via approximate SLRC), approximate gcd computation of polynomials with noisy coefficients (via SLRA), and time-series denoising (via SLRA).
\end{enumerate}
The paper is structured as follows:
\begin{enumerate}
    \item In \cref{sec:Introduction}, we introduce the necessary notations and motivate the problem.
    \item In \cref{sec:Characterization}, we characterize the critical points of the problem.
    \item In \cref{sec:LocalMinimumAlgorithm}, we propose an algorithm for finding a ``good'' local minimum of the problem.
    \item In \cref{sec:DistributionallyRobustFormulation}, we explore a distributionally robust version of the problem.
    \item In \cref{sec:NumericalResults}, we demonstrate some numerical results.
\end{enumerate}

\begin{notation}
    For \(n\in\N\), we define \(\irng{n}\defeq\Z\cap[1, n]\).
\end{notation}
Let \(m, n, p, N\in\N\), \(X\in\R^{p\times N}\), \(Y\in\R^{m\times N}\), and \(\P_{m, n, p}\defeq \R^{m\times n}\times \R^{n\times p}\times\R^{m\times p}\). For \(n_c\in\N\), let \(h:\P_{m, n, p}\to\R^{n_c}\) be defined for all \((U, V, Z)\in\P_{m, n, p}\) by
\begin{equation*}
    h_j(U, V, Z) = \innerproduct{A_{1j}}{U} + \innerproduct{A_{2j}}{V} + \innerproduct{A_{3j}}{Z} + b_j
\end{equation*}
where for each \(j\in\irng{n_c}\), \((A_{1j}, A_{2j}, A_{3j})\in\P_{m, n, p}\) and \(b_{j}\in\R\). We concern ourselves with the following optimization problem:
\begin{equation}
    \label{eqn:TheGeneralOptimizationProblem}
    \begin{aligned}
        \min_{(U, V, Z)\in\P_{m, n, p}}\quad&\quad \frac{1}{2}\norm[\Fr]{UVX - Y}^2\\
        \mathrm{s.t.}\quad\qquad&\quad Z = UV,\\
        &\quad h(U, V, Z) = 0.
    \end{aligned}
\end{equation}
The optimization problem in (\ref{eqn:TheGeneralOptimizationProblem}) is a general form of constrained matrix factorization. By modifying the problem, one can cast the optimization problem in (\ref{eqn:TheGeneralOptimizationProblem}) into several useful forms.

\subsection{Structured Low-Rank Approximation and Factorization}
Suppose that given a matrix \(Y\in\R^{m\times p}\), we wish to approximate \(Y\) by a matrix \(Z\in\R^{m\times p}\) such that \(\rank(Z) = r\) for some \(r\in\N\) and \(Z\) satisfies some affine constraints. Then the optimization problem is given as follows.
\begin{equation}
    \label{eqn:StructuredLowRankApproximation}
    \begin{aligned}
        \min_{(U, V, Z)\in\P_{m, r, p}}\quad&\quad \frac{1}{2}\norm[\Fr]{UV - Y}^2\\
        \mathrm{s.t.}\quad\qquad&\quad Z = UV,\\
        &\quad h(Z) = 0,
    \end{aligned}
\end{equation}
where \(h:\R^{m\times p}\to\R^{n_c}\) such that \(h_j\) is affine for each \(j\in\irng{n_c}\). Clearly, this is a special case of problem (\ref{eqn:TheGeneralOptimizationProblem}). This is a structured low-rank approximation problem has a wide range of applications.

If additional structure (given by affine functions) on \(U\) and \(V\) are required, this is a structured low-rank factorization problem and is given as follows.
\begin{equation}
    \label{eqn:StructuredLowRankFactorization}
    \begin{aligned}
        \min_{(U, V, Z)\in\P_{m, r, p}}\quad&\quad \frac{1}{2}\norm[\Fr]{UV - Y}^2\\
        \mathrm{s.t.}\quad\qquad&\quad Z = UV,\\
        &\quad h_1(U) = 0,\\
        &\quad h_2(V) = 0,\\
        &\quad h_3(Z) = 0,
    \end{aligned}
\end{equation}
where \(n_1, n_2, n_3\in\N\) are the number of constraints on \(U\), \(V\), and \(Z\) respectively and \(h_1:\R^{m\times n}\to\R^{n_1}\), \(h_2:\R^{n\times p}\to\R^{n_2}\), and \(h_3:\R^{m\times p}\to\R^{n_3}\) are functions with each entry being an affine function. This is also a special case of problem (\ref{eqn:StructuredLowRankApproximation}).

\subsection{Structured Low-Rank Matrix Completion}
Suppose that given a matrix \(F\in\R^{m\times p}\) with some unknown entries, we want to complete \(F\) by another matrix \(Z\in\R^{m\times p}\) while ensuring that \(\rank(Z)\leq r\) with \(r\in\N\) and keeping the known entries fixed. Let \(n_c\in\N\) and \(h:\R^{m\times p}\to\R^{n_c}\) be such that for each \(j\in\irng{n_c}\), \(h_j\) is an affine function. Suppose \(S\subseteq\irng{m}\times \irng{p}\) is the set of all known indices. So we want to solve the following optimization problem.
\begin{equation}
    \label{eqn:MatrixCompletion1}
    \begin{aligned}
        \min_{Z\in\R^{m\times p}}\quad&\quad \frac{1}{2}\norm[\Fr]{Z}^2\\
        \mathrm{s.t.}\qquad&\quad \rank(Z)\leq r,\\
        &\quad Z_{ij} = F_{ij},\qquad\forall (i, j)\in S,\\
        &\quad h(Z) = 0.
    \end{aligned}
\end{equation}
Note that \(\rank(Z)\leq r\) if and only if for \(n\in\N\) with \(n\leq r\), there exist \(U\in\R^{m\times n}\) and \(V\in\R^{n\times p}\) such that \(Z = UV\). Using this observation, we can reformulate the optimization problem in (\ref{eqn:MatrixCompletion1}) as follows.
\begin{equation}
    \label{eqn:MatrixCompletion2}
    \begin{aligned}
        \min_{\substack{n\in\N,\\ n\leq r}}\min_{\substack{U\in\R^{m\times n},\\ V\in\R^{n\times p}}}\quad&\quad \frac{1}{2}\norm[\Fr]{UV}^2\\
        \mathrm{s.t.}\quad\qquad&\quad Z = UV,\\
        &\quad Z_{ij} = F_{ij},\qquad\forall (i, j)\in S,\\
        &\quad h(Z) = 0.
    \end{aligned}
\end{equation}
Note that the inner optimization problem is of the form \cref{eqn:TheGeneralOptimizationProblem} with \(p = N\), \(X = I_p\), and the affine constraints \(h_{ij}(Z) = Z_{ij} - F_{ij} = \innerproduct{E_{ij}}{Z} - F_{ij}\) for all \((i, j)\in S\), where \(E_{ij}\) is the basis matrix of all zeros except a one at position \((i, j)\). Once for a fixed \(n\in\N\), one can solve the above optimization problem, one can simply repeat the solution for all \(n\in\irng{r}\) to obtain the optimum of the problem (\ref{eqn:MatrixCompletion2}). 


\section{Characterization of the Critical Points}
\label{sec:Characterization}
Given the importance of the optimization problem in (\ref{eqn:TheGeneralOptimizationProblem}), we wish to analyze its critical points. Note that the Lagrangian is given by
\begin{equation}
    \L(U, V, Z, M, \mu) = \frac{1}{2}\norm[\Fr]{UVX - Y}^2 + \innerproduct{M}{UV - Z} + \mu\trans h(U, V, Z),
\end{equation}
where \(M\in\R^{m\times p}\) and \(\mu\in\R^{n_c}\) are the Lagrange multipliers. Define \(\Sigma_{XX}\defeq XX\trans\), and \(\Sigma_{XY}\defeq XY\trans\), \(\Sigma_{YX}\defeq YX\trans\), and \(\Sigma_{YY}\defeq YY\trans\). Due to the KKT necessary conditions, any critical point of (\ref{eqn:TheGeneralOptimizationProblem}) must satisfy the following conditions.

    

\begin{lemma}
\label{lemma:CriticalPoints}
Let \(\pt = (U\opt, V\opt, Z\opt, M\opt, \mu\opt)\) be a critical point of (\ref{eqn:TheGeneralOptimizationProblem}). Then we have
\begin{subequations}
    \begin{align}
        &(U\opt V\opt \Sigma_{XX} - \Sigma_{YX})V\optTrans + M\opt V\optTrans + \sum_{j = 1}^{n_c}{\mu\opt_{j} A\opt_{1j}} = 0,\\
        &U\optTrans (U\opt V\opt \Sigma_{XX} - \Sigma_{YX}) + U\optTrans M\opt + \sum_{j = 1}^{n_c}{\mu\opt_{j}A\opt_{2j}} = 0,\\
        &-M\opt + \sum_{j = 1}^{n_c}{\mu\opt_{j}A\opt_{3j}} = 0,\\
        &Z\opt = U\opt V\opt,\\
        &h(U\opt, V\opt, Z\opt) = 0.
    \end{align}
\end{subequations}
\end{lemma}
\begin{proof}
Due to the KKT necessary conditions for a critical point of a constrained optimization problem, we have
\begin{subequations}
\begin{align}
    &\frac{\partial\L}{\partial U}(\pt) = (U\opt V\opt \Sigma_{XX} - \Sigma_{YX})V\optTrans + M\opt V\optTrans + \sum_{j = 1}^{n_c}{\mu\opt_{j}\frac{\partial h_{j}}{\partial U}(U\opt, V\opt, Z\opt)} = 0,\\
    &\frac{\partial\L}{\partial V}(\pt) = U\optTrans (U\opt V\opt \Sigma_{XX} - \Sigma_{YX}) + U\optTrans M\opt + \sum_{j = 1}^{n_c}{\mu\opt_{j}\frac{\partial h_{j}}{\partial V}(U\opt, V\opt, Z\opt)} = 0,\\
    &\frac{\partial\L}{\partial Z}(\pt) = -M\opt + \sum_{j = 1}^{n_1}{\mu\opt_{j}\frac{\partial h_{j}}{\partial Z}(U\opt, V\opt, Z\opt)} = 0,\\
    &Z\opt = U\opt V\opt,\\
    &h(U\opt, V\opt, Z\opt) = 0.
\end{align}
\end{subequations}
Replacing the partial derivatives with the constant matrices \(A_{1j}\), \(A_{2j}\), and \(A_{3j}\), we obtain the lemma.
\end{proof}
It is easy to see from \cref{lemma:CriticalPoints} that, in general, the Lagrange multipliers \((M\opt, \mu\opt)\), need not be unique for a given \((U\opt, V\opt, Z\opt)\).

Now we characterize the critical points of (\ref{eqn:TheGeneralOptimizationProblem}).
\begin{notation}
    For a matrix \(A\), we denote by \(A_{:, k}\) the \(k\)\textsuperscript{th} column of \(A\) and by \(A_{k,:}\) the \(k\)\textsuperscript{th} row of \(A\).
\end{notation}
\begin{lemma}[Hessian quadrature form]
    For all ordered triples \((U\opt, V\opt, Z\opt)\in\P_{m, n, p}\), \((M\opt, \mu\opt)\in\R^{m\times N}\times \R^{n_c}\), and \(\Delta\defeq(\Delta_U, \Delta_V, \Delta_Z)\in\P_{m, n, p}\), we have
    \begin{equation}
            \label{eqn:HessianQuadratureFormBeforeConstraints}
            \nabla^2_{(U, V, Z)}\L(U\opt, V\opt, Z\opt, M\opt, \mu\opt)(\Delta, \Delta) = \innerproduct{U\opt V\opt\Sigma_{X X} - \Sigma_{Y X} + M\opt}{\Delta_U\Delta_V} + \frac{1}{2}{\norm[\Fr]{(U\opt\Delta_V + \Delta_U V\opt)X}^2}.
    \end{equation}
\end{lemma}
\begin{proof}
    Note that for all \((U\opt, V\opt, Z\opt)\in\P_{m, n, p}\), \((M\opt, \mu\opt)\in\R^{m\times p}\times \R^{n_c}\), and \(\Delta\defeq(\Delta_U, \Delta_V, \Delta_Z)\in\P_{m, n, p}\), we have
    \begin{align*}
        &\L(U\opt + \Delta_U, V\opt + \Delta_V, Z\opt + \Delta_Z, M\opt, \mu\opt)\\
        &= \frac{1}{2}\norm[\Fr]{(U\opt + \Delta_U)(V\opt + \Delta_V)X - Y}^2 + \innerproduct{M\opt}{(U\opt + \Delta_U)(V\opt + \Delta_V) - Z\opt - \Delta_Z}\\&\quad + \mu\optTrans h(U\opt + \Delta_U, V\opt + \Delta_V, Z\opt + \Delta_Z)\\
        &= \frac{1}{2}\norm[\Fr]{U\opt V\opt X - Y}^2 + \innerproduct{U\opt V\opt X - Y}{(U\opt\Delta_V + \Delta_U V\opt + \Delta_U\Delta_V)X}\\&\quad + \frac{1}{2}{\norm[\Fr]{(U\opt\Delta_V + \Delta_U V\opt + \Delta_U\Delta_V)X}^2} + \innerproduct{M\opt}{U\opt V\opt - Z\opt} + \innerproduct{M\opt}{U\opt\Delta_V + \Delta_U V\opt + \Delta_U\Delta_V - \Delta_Z}\\&\quad + \mu\optTrans h(U\opt, V\opt, Z\opt) + \innerproduct{\Delta_U}{\sum_{j = 1}^{n_c}\mu_{j}\opt\frac{\partial h_{j}}{\partial U}(U\opt, V\opt, Z\opt)} + \innerproduct{\Delta_V}{\sum_{j = 1}^{n_c}\mu_{j}\opt\frac{\partial h_{j}}{\partial V}(U\opt, V\opt, Z\opt)}\\&\quad + \innerproduct{\Delta_U}{\sum_{j = 1}^{n_c}\mu_{j}\opt\frac{\partial h_{j}}{\partial Z}(U\opt, V\opt, Z\opt)}.
    \end{align*}
    Hence we have
    \begin{equation}
        \begin{aligned}
            &\L(U\opt + \Delta_U, V\opt + \Delta_V, Z\opt + \Delta_Z, M\opt, \mu\opt) - \L(U\opt, V\opt, Z\opt, M\opt, \mu\opt)\\
            &= \innerproduct{U\opt V\opt\Sigma_{XX} - \Sigma_{YX}}{U\opt\Delta_V + \Delta_U V\opt + \Delta_U\Delta_V} + \frac{1}{2}{\norm[\Fr]{(U\opt\Delta_V + \Delta_U V\opt + \Delta_U\Delta_V)X}^2}\\&\quad + \innerproduct{M\opt}{U\opt\Delta_V + \Delta_U V\opt + \Delta_U\Delta_V - \Delta_Z} + \innerproduct{\Delta_U}{\sum_{j = 1}^{n_c}\mu_{j}\opt A_{1j}} + \innerproduct{\Delta_V}{\sum_{j = 1}^{n_c}\mu_{j}\opt A_{2j}}\\&\quad + \innerproduct{\Delta_U}{\sum_{j = 1}^{n_c}\mu_{j}\opt A_{3j}}.
        \end{aligned}
    \end{equation}
    Taking only the quadratic terms, we have that
    \begin{equation}
        \nabla^2_{(U, V, Z)}\L(U\opt, V\opt, Z\opt, M\opt, \mu\opt)(\Delta, \Delta) = \innerproduct{U\opt V\opt\Sigma_{X X} - \Sigma_{Y X} + M\opt}{\Delta_U\Delta_V} + \frac{1}{2}{\norm[\Fr]{(U\opt\Delta_V + \Delta_U V\opt)X}^2},
    \end{equation}
    as claimed.
\end{proof}
\begin{definition}
    Given the optimization problem (\ref{eqn:TheGeneralOptimizationProblem}) and a pair \((U, V)\in\R^{m\times n}\times \R^{n\times p}\), define \(K_1(U, V)\in\R^{(mp)\times (mn + np + mp)}\) to be the matrix consisting of the rows
    \begin{equation}
        \begin{aligned}
            (K_1(U, V))_{(j - 1)p + k, :} = \begin{bmatrix}
                \vectorize_{m, n}(e_jV_{:, k}\trans)\trans & \vectorize_{n, p}(U_{j,:}e_k\trans)\trans & -\vectorize_{m, p}(E_{jk})\trans
            \end{bmatrix},
        \end{aligned}
    \end{equation}
    for each \((j, k)\in\irng{m}\times \irng{p}\), and define \(K_2(U, V)\in\R^{n_c\times (mn + np + mp)}\) to be the matrix consisting of the rows
    \begin{equation}
        \begin{aligned}
            (K_2(U, V))_{j, :} = \begin{bmatrix}
                \vectorize_{m, n}(A_{1j})\trans & \vectorize_{n, p}(A_{2j})\trans & \vectorize_{m, p}(A_{3j})\trans
            \end{bmatrix},
        \end{aligned}
    \end{equation}
    for each \(j\in\irng{n_c}\). Define the matrix
    \begin{equation}
        K(U, V) = \begin{bmatrix}
            K_1(U, V)\\
            K_2(U, V)
        \end{bmatrix}\in\R^{(mn + n_c)\times (mn + np + mp)}.
    \end{equation}
\end{definition}
\begin{definition}
    \label{def:BasisEtc}
    Given the critical point \(\pt = (U\opt, V\opt, Z\opt, M\opt, \mu\opt)\), define the matrix \(W \defeq U\opt V\opt\Sigma_{XX} - \Sigma_{YX} + M\opt\). Let \(\nullity(K(U\opt, V\opt)) = q\) and \(\mathfrak{B}\) be a basis for \(\Null(K(U\opt, V\opt))\). Define the set
    \begin{equation}
        \begin{aligned}
            \mathcal{B}&\defeq\{(\vectorize_{m, n}\inv(P), \vectorize_{n, p}\inv(Q), \vectorize_{m, p}\inv(R))\in\P_{m, n, p}\;|\;(P, Q, R)\in\mathfrak{B}\}.
        \end{aligned}
    \end{equation}
Consider \(q > 0\) and let \(\mathcal{B}\defeq\{(P_j, Q_j, R_j)\}_{j = 1}^{q}\subseteq\P_{m, n, p}\) be an enumeration of \(\mathcal{B}\). Define the matrix \(H(\pt; \mathcal{B})\in\R^{q\times q}\) to be the matrix with the entry at \((j, k)\in\irng{q}\times\irng{q}\) given by
\begin{equation}
    \begin{aligned}
        &(H(U, V))_{jk} = 
        \frac12\innerproduct{W}{P_jQ_k + P_kQ_j} + \frac12\innerproduct{(U\opt Q_j + P_jV\opt)X}{(U\opt Q_k + P_kV\opt)X}.
    \end{aligned}
\end{equation}
\end{definition}
\begin{theorem}
    \label{thm:CriticalPoints}
    Let \(\pt = (U\opt, V\opt, Z\opt, M\opt, \mu\opt)\) be a critical point of (\ref{eqn:TheGeneralOptimizationProblem}). If \(K(U\opt, V\opt)\) has full column-rank, then \(\pt\) is unique for fixed \((M\opt, \mu\opt)\). Otherwise, consider the matrix \(H(\pt; \mathcal{B})\) as defined above. The critical point \(\pt\) is a
    \begin{enumerate}[(i)]
        \item strict local minimum if \(H(\pt; \mathcal{B})\succ 0\), and non-strict local minimum if \(H(\pt; \mathcal{B})\succeq 0\) and \(0\in\spec(H(\pt; \mathcal{B}))\).
        \item strict local maximum if \(H(\pt; \mathcal{B})\prec 0\), and non-strict local maximum if \(H(\pt; \mathcal{B})\preceq 0\) and \(0\in\spec(H(\pt; \mathcal{B}))\).
        \item strict saddle point if \(H(\pt; \mathcal{B})\) is indefinite.
    \end{enumerate}
\end{theorem}
\begin{proof}
    We use the KKT sufficient conditions \cite[Proposition~3.3.2]{bertsekas} for optimality. We are interested in the directions \(\Delta\) such that
    \begin{align*}
        &\innerproduct{\Delta_U}{\frac{\partial f_{jk}}{\partial U}(U\opt, V\opt, Z\opt)} + \innerproduct{\Delta_V}{\frac{\partial f_{jk}}{\partial V}(U\opt, V\opt, Z\opt)}\\&\qquad + \innerproduct{\Delta_Z}{\frac{\partial f_{jk}}{\partial Z}(U\opt, V\opt, Z\opt)} = 0,\qquad\forall (j, k)\in\irng{m}\times\irng{p},\numberthis\\
        &\innerproduct{\Delta_U}{\frac{\partial h_{j}}{\partial U}(U\opt, V\opt, Z\opt)} + \innerproduct{\Delta_V}{\frac{\partial h_{j}}{\partial V}(U\opt, V\opt, Z\opt)}\\&\qquad + \innerproduct{\Delta_Z}{\frac{\partial h_{j}}{\partial Z}(U\opt, V\opt, Z\opt)} = 0,\qquad\forall j\in\irng{n_c},\numberthis
    \end{align*}
    where \(f_{jk}\) denotes the entry at \((j, k)\in\irng{m}\times\irng{p}\) of the constraint function \(f\) defined for all \((U, V, Z)\in\P_{m, n, p}\) by \(f(U, V, Z) = UV - Z\). These can be simplified as
    \begin{subequations}
        \label{eqn:FeasibleDirections}
        \begin{align*}
            &\innerproduct{\Delta_U}{e_j V_{:, k}\optTrans} + \innerproduct{\Delta_V}{U_{j, :}\opt e_k\trans} - \innerproduct{\Delta_Z}{E_{jk}} = 0,\qquad\forall (j, k)\in\irng{m}\times\irng{p},\numberthis\\
            &\innerproduct{\Delta_U}{A_{1j}} + \innerproduct{\Delta_V}{A_{2j}} + \innerproduct{\Delta_Z}{A_{3j}} = 0, \qquad\forall j\in\irng{n_c},\numberthis
        \end{align*}
    \end{subequations}
    which is equivalent to
    \begin{equation}
        K(U\opt, V\opt)\begin{bmatrix}
            \vectorize_{m, n}(\Delta_U)\\
            \vectorize_{n, p}(\Delta_V)\\
            \vectorize_{m, p}(\Delta_Z)
        \end{bmatrix} = 0.
    \end{equation}
    If \(K(U\opt, V\opt)\) has full column-rank, then the only feasible value of \(\Delta\) is \((0, 0, 0)\in\P_{m, n, p}\). So given a fixed \((M\opt, \mu\opt)\), \(\pt\) is unique. Now suppose \(q > 0\). Let \(\mathcal{B}\) be as defined in \cref{def:BasisEtc}. Then for any \(\Delta\in\P_{m, n, p}\) satisfying \cref{eqn:FeasibleDirections}, there exists \(\alpha\in\R^q\) such that
    \begin{equation}
        \Delta = \sum_{j = 1}^{q}{\alpha_j(P_j, Q_j, R_j)}.
    \end{equation}
    So we have
    \begin{subequations}
    \begin{align*}
        &\nabla^2_{(U, V, Z)}\L(\pt)(\Delta, \Delta)\\
        &= \innerproduct{W}{\sum_{j = 1}^{q}{\alpha_jP_j}\sum_{k = 1}^{q}{\alpha_kQ_k}} + \frac12\norm[\Fr]{\left(U\opt\sum_{j = 1}^{q}{\alpha_jQ_j} + \sum_{j = 1}^{q}{\alpha_jP_j}V\opt\right)X}^2\\
        &= \sum_{j = 1}^{q}{\alpha_j^2\left[\innerproduct{W}{P_{j}Q_{j}} + \frac{1}{2}\norm[\Fr]{(U\opt Q_j + P_jV\opt)X}^2\right]}\\
        &\quad + \sum_{\substack{j, k = 1\\j\neq k}}^q\alpha_j\alpha_k\bigg[\frac12\innerproduct{W}{P_jQ_k + P_kQ_j} + \frac12\innerproduct{(U\opt Q_j + P_jV\opt)X}{(U\opt Q_k + P_kV\opt)X}\bigg]\\
        &= \alpha\trans H(U\opt, V\opt; \mathcal{B})\alpha.
    \end{align*}
    \end{subequations}
    Since \(H(U\opt, V\opt; \mathcal{B})\) is a real symmetric matrix, the claims (i), (ii), and (iii) follow from the spectral theorem.
\end{proof}

\section{Algorithm for Finding a Local Minimum}
\label{sec:LocalMinimumAlgorithm}
Note that \cref{thm:CriticalPoints} characterizes the critical points of the optimization problem. Now we wish to find such critical points: in particular, local minima. Observe that if either \(U\) or \(V\) were to be made constant in (\ref{eqn:TheGeneralOptimizationProblem}), then the optimization problem becomes convex and thus every local minimum is now global. Suppose we start with an initial guess of \((U, V)\). Then we can fix \(V\) and solve for an optimal \(U\opt\) corresponding to \(V\), then fix \(U\opt\) and solve for an optimal \(V\opt\) corresponding to \(U\opt\). We can repeat this process until ``convergence'' is obtained. Since the sub-problems for fixed \(U\) or \(V\) are convex, we can use the standard methods for solving quadratic programs or the KKT necessary conditions from \cref{lemma:CriticalPoints} to solve the sub-problems. One may note that this process is a type of block-coordinate descent. Newton-Raphson method may be used after this algorithm to ensure fast convergence to a critical point. When a critical point \((U\opt, V\opt)\) is reached, we can use \cref{lemma:CriticalPoints} to compute \(M\opt\) and \(\mu\opt\). Now we can use \cref{thm:CriticalPoints} to find the nature of the critical point. If it is a local minimum, we are done. If it is a strict saddle point, we can use the method as described in the proof of \cref{thm:CriticalPoints} to find a direction \(\Delta\in\P_{m, n, p}\) that pushes the solution out of the strict saddle point. Thus, using \cref{thm:CriticalPoints} as a guide, we can always reach a local minimum. The algorithm for finding a critical point is presented in \cref{alg:CriticalPointAlgorithm}. and the algorithm for finding a local minimum is presented in \cref{alg:LocalMinimumAlgorithm}. Define \(\P^*_{m, n, p}\defeq\R^{m\times n}\times\R^{n\times p}\) and let \(L:\P^*_{m, n, p}\to\R\) be defined for all \((U, V)\in\P^*_{m, n, p}\) by \(L(U, V)\defeq\frac12\norm[\Fr]{UVX - Y}^2\).
\begin{theorem}[Termination of \cref{alg:CriticalPointAlgorithm}]
    Define the equivalence relation \(\sim\) on \(\P^*_{m, n, p}\) such that for all \((U, V), (U', V')\in\P^*_{m, n, p}\), we have \((U, V)\sim(U', V')\) if and only if \(L(U, V) = L(U', V')\) and let \([(U, V)]\) be the equivalence class of \((U, V)\). Then the while loop in \cref{alg:CriticalPointAlgorithm} converges in \({\P^*_{m, n, p}/\sim}\defeq \left\{[(U, V)]\;|\;(U, V)\in\P^*_{m, n, p}\right\}\).
\end{theorem}
\begin{proof}
    It is straight-forward to see that \(\sim\) is an equivalence relation on \(\P^*_{m, n, p}\). Let \(\{U^t, V^t\}_{t\in\NZ}\) denote the sequence of iterates produced by the while loop of \cref{alg:CriticalPointAlgorithm}. In case the algorithm terminates after \(T\in\N\) number of steps, we define \((U^t, V^t)\defeq (U^T, V^T)\) for \(t\in\N\cap(T, \infty)\). Define the sequence \(L^t\defeq\norm[\Fr]{U^tV^tX - Y}^2\) for \(t\in\NZ\). Since each sub-problem is a minimization problem, \(\{L^t\}_{t\in\NZ}\) must be a monotonically decreasing sequence bounded below by \(0\). Due to the monotone convergence theorem, the sequence \(\{L^t\}_{t\in\NZ}\) converges to its infimum. Since it is a convergent sequence, it also Cauchy. Hence there exists \(T\in\N\) such that \(L^{T - 1} - L^T < \varepsilon\). This terminates the algorithm. If the sequence is eventually constant, there exists \(T\in\N\) such that \((U^t, V^t)\in[(U^T, V^T)]\) for all \(t\in\N\cap(T, \infty]\). This completes the proof.
\end{proof}
\begin{theorem}[Termination of \cref{alg:LocalMinimumAlgorithm}]
    \cref{alg:LocalMinimumAlgorithm} terminates after a finite number of iterations.
\end{theorem}
\begin{proof}
Consider the sequence of iterates \(\{(U^t, V^t)\}_{t\in\NZ}\) produced by the algorithm. In case the algorithm terminates after \(T\in\N\) number of steps, we define \((U^t, V^t)\defeq (U^T, V^T)\) for \(t\in\N\cap(T, \infty)\). Define the sequence \(L^t\defeq\norm[\Fr]{U^tV^tX - Y}^2\) for \(t\in\NZ\). Consider the terminating conditions. If the algorithm ever encounters lines~\ref{alg:line:UniqueCriticalPoint},~\ref{alg:line:LocalMinimum} or~\ref{alg:line:NonStrictLocalMinimum}, it would terminate. So suppose it does not encounter the lines~\ref{alg:line:UniqueCriticalPoint},~\ref{alg:line:LocalMinimum} or~\ref{alg:line:NonStrictLocalMinimum}. The only other way it could terminate is through line~\ref{alg:line:WithinTolerance}. We prove that it must reach line~\ref{alg:line:WithinTolerance} eventually. Note that if it does not terminate, then it must go through line~\ref{alg:line:Descent} at each step of the iteration. Since \((\alpha\Delta_U, \alpha\Delta_V)\) in line~\ref{alg:line:Descent} is a strict descent direction, the sequence \(\{L^t\}_{t\in\NZ}\) is a strictly monotonically decreasing sequence which is bounded below by \(0\). Due to the monotone convergence theorem, as \(t\to\infty\), \(L^t\to\inf\{L^t\}_{t\in\NZ}\geq 0\) which exists as the sequence is bounded below. Since \(\{L^t\}_{t\in\NZ}\) is convergent, it must be Cauchy. Hence there exists \(T\in\N\) such that \(L^{T - 1} - L^{T} < \varepsilon\). This invokes line~\ref{alg:line:WithinTolerance} which terminates the algorithm. This completes the proof.
\end{proof}
\vfil

\begin{breakablealgorithm}
    \captionof{algorithm}{Block-coordinate descent for finding a critical point of (\ref{eqn:TheGeneralOptimizationProblem})}
    \label{alg:CriticalPointAlgorithm}
    \begin{algorithmic}[1]
        \Require{\(X\in\R^{p\times N}\), \(Y\in\R^{m\times N}\), affine function \(h:\P_{m, n, p}\to\R^{n_c}\), initial guess \((U, V)\in\P^*_{m, n, p}\) for \(m, n, p, N, n_c\in\N\), and tolerance \(\varepsilon > 0\).}
        \Ensure{A critical point \((U, V)\in\P^*_{m, n, p}\) for the optimization problem in (\ref{eqn:TheGeneralOptimizationProblem}) within a tolerance of \(\varepsilon > 0\).}
        \Procedure{Find-critical-point}{$X$, $Y$, $h$, $U$, $V$, $\varepsilon$}
        \State \(\delta\gets\infty\)
        \While {\(\delta \geq \varepsilon\)}
        \State With \(U\) fixed, compute an optimal \(V\opt\) from (\ref{eqn:TheGeneralOptimizationProblem})
        \State \(\delta\gets\norm[\mathrm{max}]{V - V\opt}\)
        \State\(V\gets V\opt\)
        \State With \(V\) fixed, compute an optimal \(U\opt\) from (\ref{eqn:TheGeneralOptimizationProblem})
        \State \(\delta\gets\max\left(\left\{\norm[\mathrm{max}]{U - U\opt}, \delta\right\}\right)\)
        \State \(U\gets U\opt\)
        \EndWhile
        \State\Return \((U, V)\)
        \EndProcedure
    \end{algorithmic}
\end{breakablealgorithm}
\vspace*{\baselineskip}
\begin{breakablealgorithm}
    \captionof{algorithm}{Finding a local minimum of (\ref{eqn:TheGeneralOptimizationProblem})}
    \label{alg:LocalMinimumAlgorithm}
    \begin{algorithmic}[1]
        \Require{\(X\in\R^{p\times N}\), \(Y\in\R^{m\times N}\), affine function \(h:\P_{m, n, p}\to\R^{n_c}\) for \(m, n, p, N, n_c\in\N\), and tolerance \(\varepsilon > 0\).}
        \Ensure{A local minimum \((U, V)\in\P^*_{m, n, p}\) for the optimization problem (\ref{eqn:TheGeneralOptimizationProblem}) within a tolerance of \(\varepsilon > 0\).}
        \Procedure{Find-local-minimum}{$X$, $Y$, $h$, $\varepsilon$}
        \State \(\Sigma_{XX}\gets XX\trans\)
        \State \(\Sigma_{YX}\gets YX\trans\)
        \State Initialize \(U\in\R^{m\times n}\), \(V\in\R^{n\times p}\)
        \State \(s\gets 0\)
        \While {\(s\neq 1\)}
        \State \((U\opt, V\opt)\gets \Call{Find-critical-point}{X, Y, h, U, V, \varepsilon}\)
        \State Form the matrix \(K(U\opt, V\opt)\) and find the basis \(\mathfrak{B}\) of \(\Null(K(U\opt, V\opt))\).
        \If {\(\mathfrak{B} = \varnothing\)}
        \State \(s\gets 1\)\label{alg:line:UniqueCriticalPoint}
        \Else
        \State Solve for \(M\opt\) and \(\mu\opt\) using \cref{lemma:CriticalPoints}.
        \State \(\pt\gets (U\opt, V\opt, U\opt V\opt, M\opt, \mu\opt)\)
        \State Form the matrix \(H(\pt; \mathcal{B})\)
        \If {\(H(\pt; \mathcal{B})\succ 0\)}
        \State \(s\gets 1\)\label{alg:line:LocalMinimum}\Comment{Strict local minimum}
        \ElsIf {\(H(\pt; \mathcal{B})\) has a negative eigenvalue \(\lambda\)}
        \State Let the ordered triple \((\vectorize_{m, n}(\Delta_U), \vectorize_{n, p}(\Delta_V), \vectorize_{m, n}(\Delta_Z))\) be an eigenvector
        \Statex \qquad\qquad\qquad\qquad corresponding to \(\lambda\).
        \Else
        \State \(s\gets 1\)\Comment{Non-strict local minimum}\label{alg:line:NonStrictLocalMinimum}
        \EndIf
        \If {\(s\neq 1\)}
        \State Choose \(\alpha > 0\) using line search until
        \State \hskip\algorithmicindent \(\norm[\Fr]{(U\opt + \alpha\Delta_U)(V\opt + \alpha\Delta_V)X - Y} < \norm[\Fr]{U\opt V\opt X - Y}\).
        \State \((U, V)\gets (U\opt + \alpha\Delta_U, V\opt + \alpha\Delta_V)\)\label{alg:line:Descent}
        \If {\(\norm[\Fr]{U\opt V\opt X - Y}^2 - \norm[\Fr]{UVX - Y}^2 < \varepsilon\)}
        \State \((U\opt, V\opt)\gets \Call{Find-critical-point}{X, Y, h, U, V, \varepsilon}\)
        \State \(s\gets 1\)\label{alg:line:WithinTolerance}
        \EndIf
        \EndIf
        \EndIf
        \EndWhile
        \State\Return \((U\opt, V\opt)\)
        \EndProcedure
    \end{algorithmic}
\end{breakablealgorithm}

\section{Distributionally Robust Formulation}
\label{sec:DistributionallyRobustFormulation}
Suppose the pair \(\xi \defeq [\begin{matrix}X\trans & Y\trans\end{matrix}]\in\R^{N\times(m + p)}\) in problem~\ref{eqn:TheGeneralOptimizationProblem} is random and it comes from a probability space \((\Omega, \mathscr{B}, P)\). We shall further assume that \(\Omega\) is Polish and denote the induced probability measure by \(\Pr = \xi_\#P\) which is unknown. Denote \(\Xi\defeq\supp\Pr\).

Suppose we have \(M\in\N\) samples of \(\xi\) which we denote by \(\{\xi^j\}_{j = 1}^M\). We want to transform problem~\ref{eqn:TheGeneralOptimizationProblem} into a distributionally robust formulation. Since \(\xi\) does not appear in the constraints, we only look at the objective. For \(r\in(1, \infty)\), we consider all distributions in the ambiguity set given by the Wasserstein \(r\)-ball \(\bar B(\widehat\Pr_M; \varepsilon)\) of radius \(\varepsilon > 0\) around the empirical distribution \(\widehat\Pr_M\defeq\frac1M\sum_{j = 1}^{M}\delta_{\xi^j}\).

We minimize the worst-case value of the objective with respect to the distributions in the ambiguity set and obtain the optimization problem
\begin{equation}
    \label{eqn:RobustReformulation}
    \begin{aligned}
        \min_{(U, V, Z)\in\P_{m, n, p}}\quad&\quad \sup_{\Q\in\bar B(\widehat\Pr_M; \varepsilon)}\E_\Q\norm[\Fr]{UVX - Y}\\
        \mathrm{s.t.}\quad\qquad&\quad Z = UV,\\
        &\quad h(U, V, Z) = 0.
    \end{aligned}
\end{equation}
Note that \(\norm[\Fr]{ZX - Y} = \norm[\Fr]{X\trans Z\trans - Y\trans} = \norm{\tilde{\xi} \tilde{\eta}}\), where \(\eta\defeq [\begin{matrix}Z & -I_m\end{matrix}]\trans\), and \(\tilde\eta\defeq \Vect_{(m + p), m}(\eta)\), and \(\tilde\xi\in\R^{mN\times(m^2 + mp)}\) is such that \(\Vect_{N, m}(\xi\eta) = \tilde{\xi}\tilde{\eta}\). 
From standard results \cite{esfahani-kuhn}, we have
\begin{equation}
        \sup_{\Q\in\bar B(\widehat\Pr_M; \varepsilon)}\E_\Q\norm{\tilde{\xi}\tilde{\eta}} = \inf_{\lambda \geq 0}\left\{\lambda\varepsilon + \frac1M\sum_{j = 1}^{M}\sup_{\xi\in\Xi}\left\{\norm{\tilde\xi\tilde\eta} - \lambda\norm[r]{\tilde\xi - \tilde\xi^j}\right\}\right\}.
\end{equation}

Since the map \(x\mapsto\norm{x}\) is convex and \(1\)-Lipschitz, it is well-known \cite{coulson-2021} that
\begin{equation}
    \sup_{\xi\in\Xi}\left\{\norm{\tilde\xi\tilde\eta} - \lambda\norm[r]{\tilde\xi - \tilde\xi^j}\right\}\leq\begin{cases}
        \norm{\tilde\xi^j\tilde\eta},\quad&\text{if } \norm[q]{\tilde{\eta}}\leq\lambda\\
        \infty,\quad&\text{otherwise,}
    \end{cases}
\end{equation}
where \(q\) is such that \(\frac1r + \frac1q = 1\). Note that if \(\Xi = \R^{N\times(m + p)}\), then the above inequality holds with equality. So we have the relaxed problem with the affine map \(Z\mapsto\tilde\eta\)
\begin{equation}
    \label{eqn:RelaxedProblem}
    \begin{aligned}
        \min_{\substack{(U, V, Z)\in\P_{m, n, p},\\\lambda\in\R}}\quad&\quad \lambda\varepsilon + \frac1M\sum_{j = 1}^{M}\norm{\tilde\xi^j\tilde\eta}\\
        \mathrm{s.t.}\quad\qquad&\quad Z = UV,\\
        &\quad h(U, V, Z) = 0,\\
        &\quad\lambda\geq 0,\\
        &\quad\norm[r,*]{\tilde\eta}\leq\lambda,
    \end{aligned}
\end{equation}
which is equivalent to the original problem~\ref{eqn:RobustReformulation} if \(\Pr\) has full support. Define the affine function \(\varphi:\R^{m\times p}\to\R^{m^2 + mp}\) given for all \(Z\in\R^{m\times p}\) by \(\varphi(Z) = \Vect_{(m + p), m}\left([\begin{matrix}Z & -I_m\end{matrix}]\trans\right)\). Then problem~\ref{eqn:RelaxedProblem} is equivalent to the optimization problem
\begin{equation}
    \label{eqn:RelaxedProblemEquivalent}
    \begin{aligned}
        \min_{(U, V, Z)\in\P_{m, n, p}}\quad&\quad \frac1M\sum_{j = 1}^{M}\norm{\tilde\xi^j\varphi(Z)} + \varepsilon\norm[r,*]{\varphi(Z)}\\
        \mathrm{s.t.}\quad\qquad&\quad Z = UV,\\
        &\quad h(U, V, Z) = 0.
    \end{aligned}
\end{equation}
Since problem~\ref{eqn:RelaxedProblemEquivalent} is convex for fixed \(U\) or \(V\), \cref{alg:CriticalPointAlgorithm} works as expected. However, \cref{alg:LocalMinimumAlgorithm} no longer works readily. Furthermore, since problem~\ref{eqn:RelaxedProblemEquivalent} no longer follows the formulation of problem~\ref{eqn:TheGeneralOptimizationProblem}, we may relax problem~\ref{eqn:RelaxedProblemEquivalent} in the case of \(r = 2\) by squaring each of the individual terms in the objective, and noting that \(\norm{\varphi(Z)}^2 = \norm[\Fr]{Z}^2 + \norm[\Fr]{-I_m}^2 = \norm[\Fr]{Z}^2 + m\), which results in the following optimization problem:
\begin{equation}
    \label{eqn:MoreRelaxedProblem}
    \begin{aligned}
        \min_{(U, V, Z)\in\P_{m, n, p}}\quad&\quad \frac{1}{2M}\sum_{j = 1}^{M}\norm[\Fr]{ZX^j - Y^j}^2 + \frac{\tilde\varepsilon}{2}\norm[\Fr]{Z}^2\\
        \mathrm{s.t.}\quad\qquad&\quad Z = UV,\\
        &\quad h(U, V, Z) = 0,
    \end{aligned}
\end{equation}
which follows the formulation of problem~\ref{eqn:TheGeneralOptimizationProblem} with \(\tilde\varepsilon > 0\) correspondingly being the new hyperparameter instead of \(\varepsilon\). Thus, we have arrived at the regularized multiple-sample version of problem~\ref{eqn:TheGeneralOptimizationProblem}.

\section{Numerical Results}
\label{sec:NumericalResults}
\subsection{The Initial Values}
In some cases, choice of a good initial value can speed up convergence drastically. To find good initial values, we start with either a randomly initialized value or matrices with all diagonals equal to \(1\), follow \cref{alg:CriticalPointAlgorithm}, and after each step, we add a small noise sampled i.i.d. from uniform distributions to the iterates. We decay the standard deviation of the sampling distribution by a constant factor in \((0, 1)\) after each iteration. We stop when the standard deviation of the noise goes below a certain threshold. This selects a good initial value by virtue of the nature of \cref{alg:CriticalPointAlgorithm}.

Another way of finding good initial conditions is to run a fast, simple algorithm such as Cadzow for a few iterations to obtain an approximate \(\widehat Z\) and then factorize this \(\hat Z\) into two matrices \(\widehat U\) and \(\widehat V\) of the required shape (via e.g. SVD). This ``warm-start'' dramatically speeds up convergence, especially for low-noise situations.

\subsection{SLRA Sylvester Matrices}
Suppose we have an LTI system given by the transfer function \(G : \C\to\C^*\) given for all \(s\in\C\) by \(G(s) = \frac{p(s)}{q(s)}\) where \(p, q\in\R[s]\). It is of great interest to know if \(p\) or \(q\) share a zero, i.e., \(\deg(\gcd(p, q))\). In reality, due to random noise, this degree will almost surely be \(1\). We want to approximate the Sylvester matrix generated by \(p\) and \(q\) by another low-rank Sylvester matrix. We take two polynomials given for all \(s\in\C\) by
\begin{subequations}
    \begin{align}
        p(s) &= 5 - 6s + s^2 = (5 - s)(1 - s),\\
        q(s) &= 5.72 - 6.3s  + s^2 = (5.2 - s)(1.1 - s).
    \end{align}
\end{subequations}
We interpret the polynomial \(q\) as a perturbed version of \(p\). Clearly, the perturbation changes the zeros of \(q\) from that of \(p\), i.e., \(\gcd(p, q) = 1\). We want to approximate the polynomials \(p\) and \(q\) by \(\hat p\in\R[s]\), and \(\hat q\in\R[s]\) respectively such that \(\deg(\gcd(\hat p, \hat q)) = 1\). To do this, we form the Sylvester matrix
\begin{equation*}
    S = \begin{bmatrix}
        1.00 & -6.00 &  5.00 & 0.00\\
        0.00 &  1.00 & -6.00 & 5.00\\
        1.00 & -6.30 &  5.72 & 0.00\\
        0.00 &  1.00 & -6.30 & 5.72
    \end{bmatrix}
\end{equation*}
from \(p\) and \(q\), and perform a structured low-rank approximation with rank equal to \(4 - 1 = 3\). This produces the approximation
\begin{equation}
    \widehat{S} =
    \begin{bmatrix}
        0.98499759 & -6.00294193 & 4.99942333 & 0\\
        0 &  0.98499759 & -6.00294193 & 4.99942333\\
        1.01494403 & -6.2970688 & 5.7205752 & 0\\
        0 & 1.01494403 & -6.2970688 & 5.7205752
    \end{bmatrix},
\end{equation}
and correspondingly
\begin{equation}
        \begin{aligned}
            \hat p(s) &\approx (5.09895949 - s)(0.99541271 - s),\\
            \hat q(s) &\approx (5.09895949 - s)(1.1053913 - s),
        \end{aligned}
\end{equation}
with \(\deg(\gcd(\hat{p}, \hat{q})) = 1\) as expected. Note that \(\norm{\hat p - p}^2 + \norm{\hat q - q}^2 = \frac12\norm[\Fr]{\widehat{S} - S}^2 \approx 0.0004663\) and \(\rank(\widehat{S}) = 3\). This matches the results obtained in \cite{zhi-approximate-gcd} in the case when the leading coefficient is allowed to be perturbed.

\subsection{SLRA of Hankel Matrices}
\label{subsec:SLRAHankel}
We take the time-series \(\{y_t\}_{t\in\irng{19}}\), given for all \(t\in\irng{19}\) by \(y_t = e^{\lambda t}\sin(\omega t)\), where \(\lambda = 0.05\) and \(\omega = 0.4\pi\), and add ``noise'' terms \(\{\varepsilon_t\}_{t\in\irng{19}}\) to get the time-series \(\{\tilde{y}_t\}_{t\in\irng{19}}\) given for all \(t\in\irng{19}\) by \(\tilde{y}_t\defeq y_t + \varepsilon_t\). We try to reconstruct \(\{y_t\}_{t\in\irng{19}}\) from \(\{\tilde{y}_t\}_{t\in\irng{19}}\). We form a Hankel matrix \(H\in\R^{10\times 10}\) from \(\{\tilde{y}_t\}_{t\in\irng{19}}\) and find a low-rank approximation with rank \(2\). We denote the reconstructed time-series by \(\{\hat{y}_t\}_{t\in\irng{19}}\).

\begin{figure}[H]
    \centering
    \begin{tikzpicture}
    \begin{axis}[
        boxplot/every box/.style = {thick, solid, red},
        boxplot/every median/.style = {thick, solid, blue},
        boxplot/every whisker/.style = {thick, solid, blue},
        boxplot/draw direction = y,
        grid = both,
        grid style = {dashed},
        minor grid style = {dashed},
        width = 8.5cm,
        height = 6cm,
        xlabel = {\(\sigma\)},
        ylabel = {\(\norm{\hat{y} - y}\)},
        xtick = {0.3, 0.6, 0.9}
    ]
        \addplot+ [boxplot prepared = {
        	draw position = 0.30000,
        	lower whisker = 0.31775, lower quartile = 0.52249,
        	median = 0.70467,
        	upper quartile = 0.79371, upper whisker = 0.89036,
        	box extend = 0.02500
        	},
        ] coordinates {};
        \addplot+ [boxplot prepared = {
        	draw position = 0.60000,
        	lower whisker = 1.00706, lower quartile = 1.07410,
        	median = 1.17060,
        	upper quartile = 1.76246, upper whisker = 2.06103,
        	box extend = 0.02500
        	},
        ] coordinates {};
        \addplot+ [boxplot prepared = {
        	draw position = 0.90000,
        	lower whisker = 1.31921, lower quartile = 1.50288,
        	median = 1.91345,
        	upper quartile = 2.13330, upper whisker = 2.59624,
        	box extend = 0.02500
        	},
        ] coordinates {};
    \end{axis}
    \end{tikzpicture}
    \caption{Error between the true and the reconstructed time series.}
    \label{fig:HankelMultipleRuns}
\end{figure}
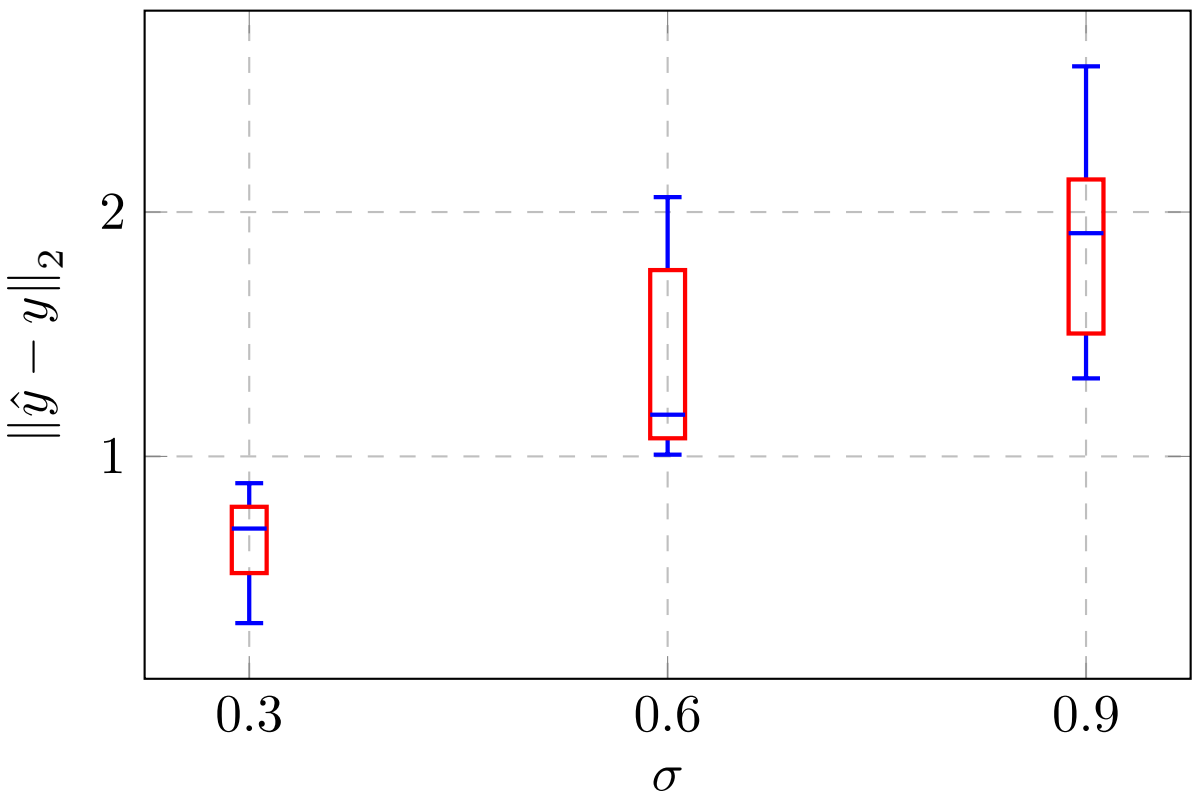

We first choose \(\varepsilon_t = (-1)^t c\) for all \(t\in\irng{19}\), and choose \(c\in\{0.2, 0.4, 0.6, 0.8\}\), and obtain the approximations for each of them. The reconstruction errors are given in \cref{tab:HankelDeterministicNoise} and the resulting time-series are shown in \cref{fig:HankelDeterministicNoise}. Next, we sample \(\{\varepsilon_t\}_{t\in\irng{19}}\) i.i.d. from \(\mathcal{N}(0, \sigma^2)\) for \(\sigma\in\{0.3, 0.6, 0.9\}\), each \(10\) times and compute the reconstructed time-series. The resulting reconstruction errors are given in \cref{fig:HankelMultipleRuns}.

\begin{table}[H]
    \centering
    \captionof{table}{Reconstruction errors for different values of \(c\).}
    \begin{tabular}{c c c}
        \toprule
        \(c\) & \(\norm{\tilde{y} - y}\) & \(\norm{\hat{y} - y}\)\\
        \midrule
        0.1 & 0.43589 & 0.03722\\
        0.2 & 0.87178 & 0.07446\\
        0.3 & 1.30767 & 0.11246\\
        0.4 & 1.74356 & 0.15191\\
        0.5 & 2.17945 & 0.19543\\
        0.6 & 2.61534 & 0.25148\\
        0.7 & 3.05123 & 0.36764\\
        0.8 & 3.48712 & 1.10242\\
        \bottomrule
    \end{tabular}
    \label{tab:HankelDeterministicNoise}
\end{table}

\begin{figure}[H]
    \centering
    \begin{tikzpicture}
    \begin{groupplot}[
        group style = {
            group size = 2 by 1,
            horizontal sep = 1.5cm,
            vertical sep = 1.3cm,
        },
        auto title,
    ]
    \nextgroupplot[
        axis lines = left,
        xlabel = {Time (\(t\))},
        ylabel = {Time series},
        xmajorgrids,
        ymajorgrids,
        grid style = dashed,
        width = 8cm,
        height = 6cm,
        legend style = {
            at = {(0.05, 0.975)},
            anchor = north west,
            nodes = {scale = 1.0, transform shape}
        }
    ]
        \addplot+[no markers, very thick, color = blue] table [x = t, y = y, col sep = comma] {data/Hankel_Cadzow_deterministic_c_0.6.csv};
        \addplot+[no markers, very thick, color = red] table [x = t, y = y_noisy, col sep = comma] {data/Hankel_Cadzow_deterministic_c_0.6.csv};
        \addplot+[no markers, very thick, color = green] table [x = t, y = y_rec, col sep = comma] {data/Hankel_Cadzow_deterministic_c_0.6.csv};
        \legend{\(y\), \(\tilde{y}\), \(\hat{y}\)}
    \nextgroupplot[
        axis lines = left,
        xlabel = {Time (\(t\))},
        ylabel = {Time series},
        xmajorgrids,
        ymajorgrids,
        grid style = dashed,
        width = 8cm,
        height = 6cm,
        legend style = {
            at = {(0.05, 0.975)},
            anchor = north west,
            nodes = {scale = 1.0, transform shape}
        }
    ]
        \addplot+[no markers, very thick, color = blue] table [x = t, y = y, col sep = comma] {data/Hankel_Cadzow_deterministic_c_0.8.csv};
        \addplot+[no markers, very thick, color = red] table [x = t, y = y_noisy, col sep = comma] {data/Hankel_Cadzow_deterministic_c_0.8.csv};
        \addplot+[no markers, very thick, color = green] table [x = t, y = y_rec, col sep = comma] {data/Hankel_Cadzow_deterministic_c_0.8.csv};
        \legend{\(y\), \(\tilde{y}\), \(\hat{y}\)}
    \end{groupplot}
    \end{tikzpicture}
    \captionof{figure}{The true (\(y\)), noisy (\(\tilde{y}\)), and reconstructed (\(\hat{y}\)) time-series for \(c = 0.6\) (a), and \(c = 0.8\) (b).}
    \label{fig:HankelDeterministicNoise}
\end{figure}
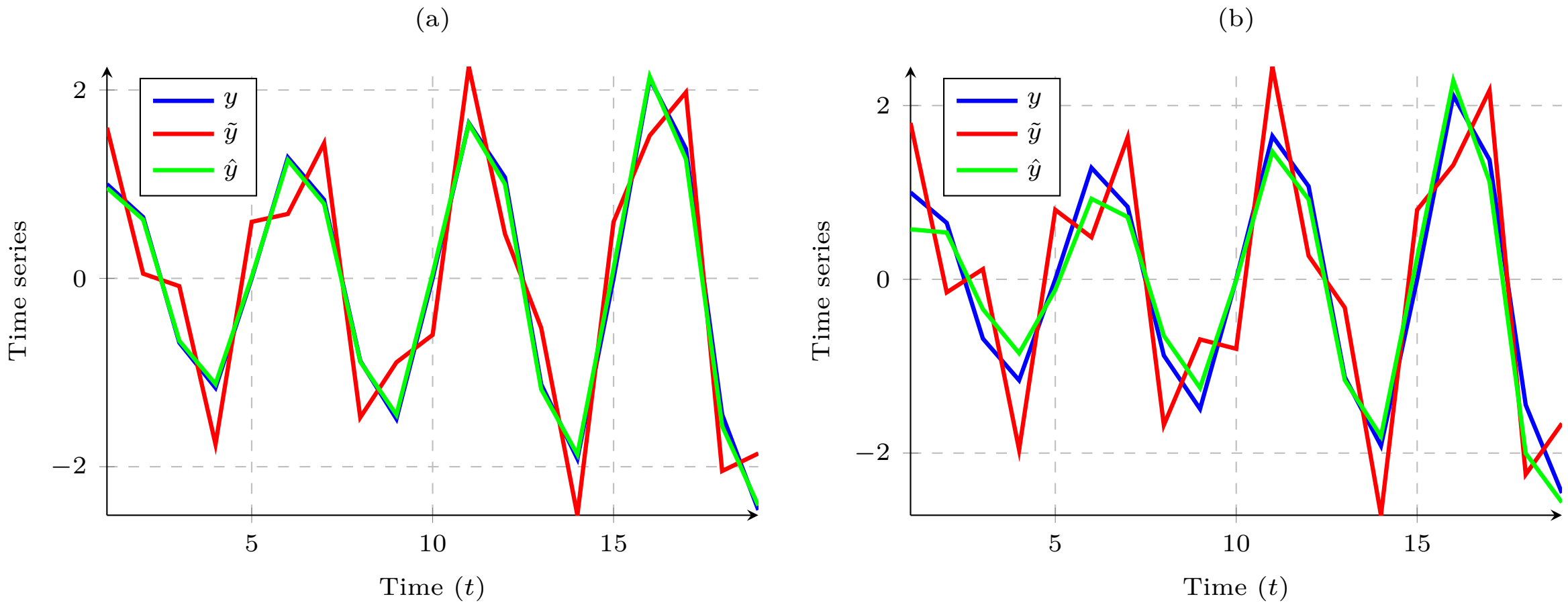

\subsection{Unstructured Matrix Completion}
We take the first \(100\times 100\) block, say \(Y\in\R^{100\times 100}\), of the traffic dataset \cite{traffic}, mask \(\SI{50}{\%}\) of the entries and perform a low-rank completion of the masked matrix with rank \(r = 8\), which captures \(\SI{99}{\%}\) of the variance of the unmasked matrix. This results in a matrix \(Z\in\R^{100\times 100}\) with an MSE of \(\frac{1}{10^4}\norm[\Fr]{Z - Y}^2\approx 5.0121\times 10^{-5}\). Performing a low-rank decomposition with NNM and SVT yield MSEs of \(8.1412\times 10^{-4}\) and \(7.7855\times 10^{-4}\) respectively. The \(90\)\textsuperscript{th} quantile values of the squared errors are approximately \(6.0502\times 10^{-5}\), \(2.5530\times 10^{-3}\), and \(2.4803\times 10^{-3}\).

\subsubsection{Asymptotic Complexity}
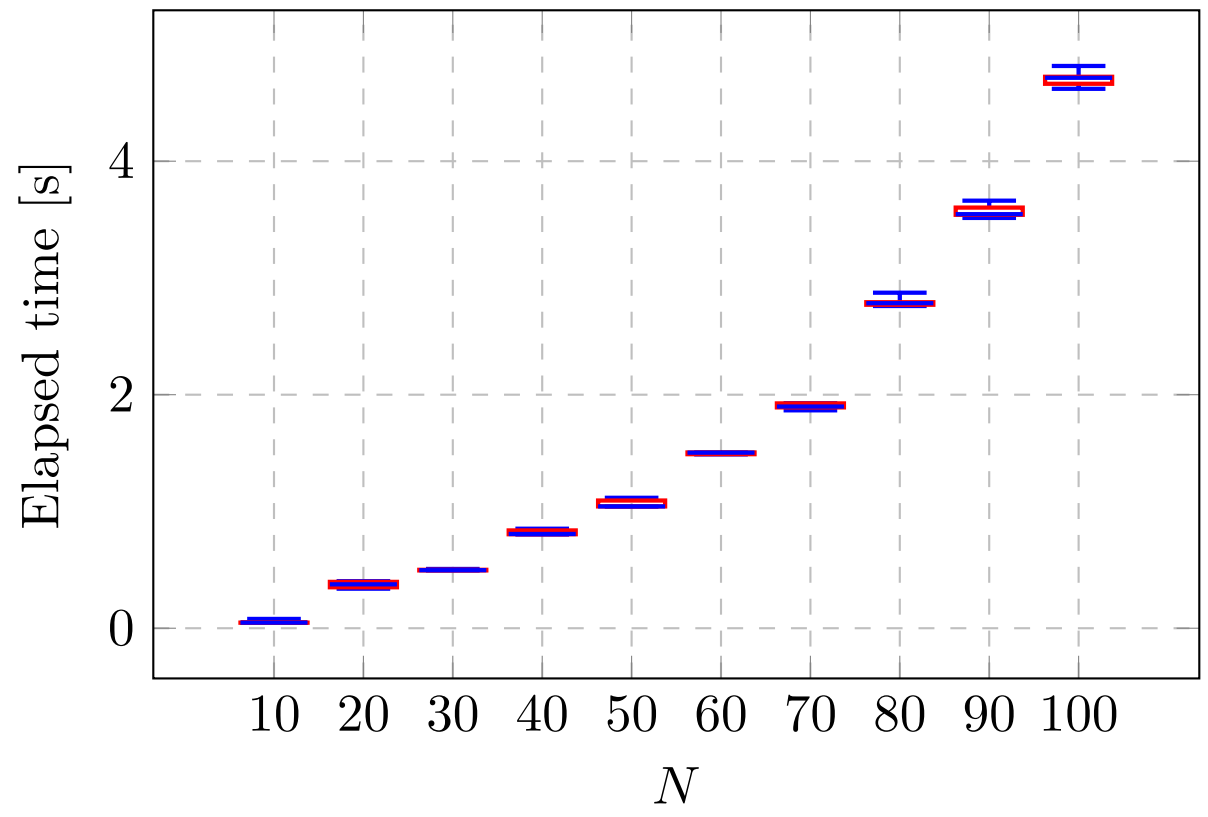
\begin{figure}[H]
    \centering
    \begin{tikzpicture}
    \begin{axis}[
        boxplot/every box/.style = {thick, solid, red},
        boxplot/every median/.style = {thick, solid, blue},
        boxplot/every whisker/.style = {thick, solid, blue},
        boxplot/draw direction = y,
        grid = both,
        grid style = {dashed},
        minor grid style = {dashed},
        width = 8.5cm,
        height = 6cm,
        xlabel = {\(N\)},
        ylabel = {Elapsed time [\(\si{\second}\)]},
        xtick = {10, 20, 30, 40, 50, 60, 70, 80, 90, 100}
    ]
        \addplot+ [boxplot prepared = {
            draw position = 10.00000,
            lower whisker = 0.04685, lower quartile = 0.04719,
            median = 0.04797,
            upper quartile = 0.04814, upper whisker = 0.07956,
            box extend = 7.5
            },
        ] coordinates {};
        \addplot+ [boxplot prepared = {
            draw position = 20.00000,
            lower whisker = 0.33850, lower quartile = 0.35072,
            median = 0.37592,
            upper quartile = 0.39704, upper whisker = 0.40287,
            box extend = 7.5
            },
        ] coordinates {};
        \addplot+ [boxplot prepared = {
            draw position = 30.00000,
            lower whisker = 0.49470, lower quartile = 0.49500,
            median = 0.49723,
            upper quartile = 0.50038, upper whisker = 0.50759,
            box extend = 7.5
            },
        ] coordinates {};
        \addplot+ [boxplot prepared = {
            draw position = 40.00000,
            lower whisker = 0.80258, lower quartile = 0.80481,
            median = 0.80704,
            upper quartile = 0.83760, upper whisker = 0.85219,
            box extend = 7.5
            },
        ] coordinates {};
        \addplot+ [boxplot prepared = {
            draw position = 50.00000,
            lower whisker = 1.04153, lower quartile = 1.04312,
            median = 1.04341,
            upper quartile = 1.09371, upper whisker = 1.11715,
            box extend = 7.5
            },
        ] coordinates {};
        \addplot+ [boxplot prepared = {
            draw position = 60.00000,
            lower whisker = 1.48879, lower quartile = 1.49044,
            median = 1.50369,
            upper quartile = 1.50546, upper whisker = 1.50546,
            box extend = 7.5
            },
        ] coordinates {};
        \addplot+ [boxplot prepared = {
            draw position = 70.00000,
            lower whisker = 1.86610, lower quartile = 1.89122,
            median = 1.89868,
            upper quartile = 1.92551, upper whisker = 1.92679,
            box extend = 7.5
            },
        ] coordinates {};
        \addplot+ [boxplot prepared = {
            draw position = 80.00000,
            lower whisker = 2.75969, lower quartile = 2.77043,
            median = 2.78375,
            upper quartile = 2.79215, upper whisker = 2.87360,
            box extend = 7.5
            },
        ] coordinates {};
        \addplot+ [boxplot prepared = {
            draw position = 90.00000,
            lower whisker = 3.51399, lower quartile = 3.54082,
            median = 3.54615,
            upper quartile = 3.60259, upper whisker = 3.66181,
            box extend = 7.5
            },
        ] coordinates {};
        \addplot+ [boxplot prepared = {
            draw position = 100.00000,
            lower whisker = 4.61981, lower quartile = 4.66290,
            median = 4.71429,
            upper quartile = 4.72301, upper whisker = 4.81551,
            box extend = 7.5
            },
        ] coordinates {};
    \end{axis}
    \end{tikzpicture}
    \caption{Time elapsed for \(10\) iterations of the algorithm for the first \(N\times N\) block of the traffic dataset.}
    \label{fig:UnstructuredComplexity}
\end{figure}

We take the \(N\times N\) block of the traffic dataset and run \(10\) iterations of the algorithm \(5\) times for each \(N\in\{10j\;|\;j\in\irng{10}\}\) and measure the elapsed times. The results are shown in \cref{fig:UnstructuredComplexity}. By a least-squares fit after taking the logarithm of the elapsed times, the elapsed time seems to grow as \(\Order\left(N^{1.8}\right)\).

\subsection{Structured Matrix Completion}
Here, we are interested in completing Hankel matrices. We take the same time-series as described in \cref{subsec:SLRAHankel} and delete some of the entries using i.i.d. \(\mathrm{Bernoulli}(p)\) trials where \(p\in(0, 1)\) is the probability of omission. We construct a Hankel matrix using this masked time-series \(\tilde{y}\) and perform SLRC with rank \(r = 2\) to get the reconstructed time-series \(\hat y\). The results are shown in \cref{fig:HankelSLRC0.7} and the errors are given in \cref{tab:HankelSLRCError}.

Alternatively, one may solve the problem
\begin{equation}
    \label{eqn:MatrixCompletionAlternative1}
    \begin{aligned}
        \min_{Z\in\R^{m\times p}}\quad&\quad \frac{1}{2}\norm[\Fr]{(Z - F)\odot \Omega}^2\\
        \mathrm{s.t.}\qquad&\quad \rank(Z)\leq r,\\
        &\quad h(Z) = 0,
    \end{aligned}
\end{equation}
where \(\Omega\in\R^{m\times n}\) such that for all \((i, j)\in\irng{m}\times\irng{n}\),
\begin{equation}
    \Omega_{ij} = \begin{cases}
        1, \quad &\text{if } (i, j)\in S\\
        0, \quad &\text{if } (i, j)\notin S.
    \end{cases}
\end{equation}
This can be recast as
\begin{equation}
    \label{eqn:MatrixCompletionAlternative2}
    \begin{aligned}
        \min_{\substack{n\in\N,\\ n\leq r}}\min_{\substack{U\in\R^{m\times n},\\ V\in\R^{n\times p}}}\quad&\quad \frac{1}{2}\norm[\Fr]{(UV - F)\odot \Omega}^2\\
        \mathrm{s.t.}\quad\qquad&\quad Z = UV,\\
        &\quad h(Z) = 0,
    \end{aligned}
\end{equation}
with minor changes in the algorithm.
\begin{remark}
    If the known entries are particularly noisy, then a crude initial guess may result in infeasible problems if we straight-up use (\ref{eqn:MatrixCompletion2}). It may be beneficial in that scenario to use (\ref{eqn:MatrixCompletionAlternative2}) to arrive at a feasible solution of (\ref{eqn:MatrixCompletion2}) and use (\ref{eqn:MatrixCompletion2}) henceforth.
\end{remark}
\begin{table}[H]
    \centering
    \captionof{table}{The error between the original and completed time-series.}
    \begin{tabular}{c c}
        \toprule
        \(p\) & \(\norm{y - \hat y}\)\\
        \midrule
        \(0.1\) & \(0.00015\)\\
        \(0.2\) & \(0.00100\)\\
        \(0.5\) & \(0.00144\)\\
        \(0.7\) & \(0.05813\)\\
        \bottomrule
    \end{tabular}
    \label{tab:HankelSLRCError}
\end{table}
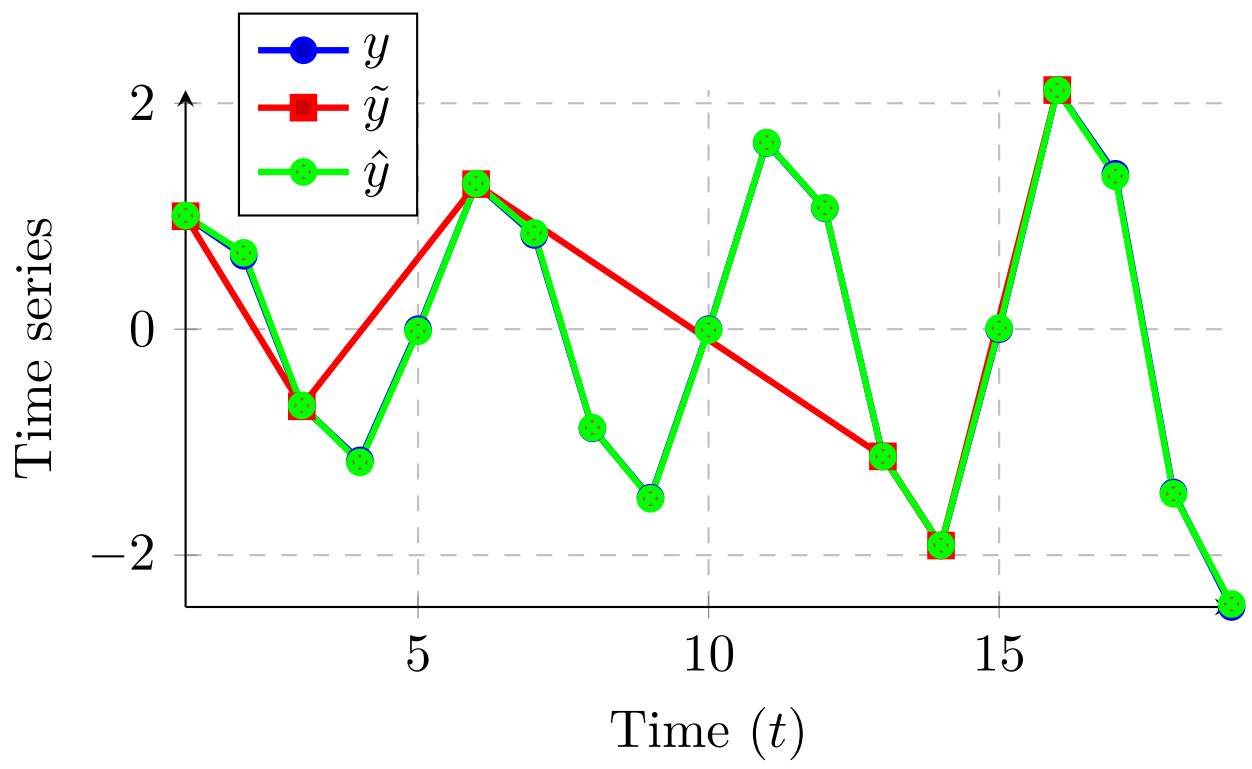
\begin{figure}[H]
    \centering
    \begin{tikzpicture}
        \begin{axis}[
            axis lines = left,
            xlabel = {Time (\(t\))},
            ylabel = {Time series},
            xmajorgrids,
            ymajorgrids,
            grid style = dashed,
            width = 8.5cm,
            height = 5cm,
            legend style = {
                at = {(0.05, 1.15)},
                anchor = north west
            }
        ]
            \addplot+[very thick, color = blue] table [x = t, y = y, col sep = comma] {data/Hankel_omit_p_0.7.csv};
            \addplot+[very thick, color = red] table [x = t, y = y_masked, col sep = comma] {data/Hankel_omit_p_0.7.csv};
            \addplot+[very thick, color = green] table [x = t, y = y_rec, col sep = comma] {data/Hankel_omit_p_0.7.csv};
            \legend{\(y\), \(\tilde{y}\), \(\hat{y}\)}
        \end{axis}
    \end{tikzpicture}
    \captionof{figure}{The true (\(y\)), masked (\(\tilde{y}\)), and reconstructed (\(\hat{y}\)) time-series for \(p = 0.7\).}
    \label{fig:HankelSLRC0.7}
\end{figure}

\subsubsection{Time-Series Forecasting Using Approximate Matrix Completion}
\begin{figure}[H]
    \centering
    \begin{tikzpicture}
        \begin{axis}[
            axis lines = left,
            xlabel = {Time (\(t\))},
            ylabel = {Time series},
            xmajorgrids,
            ymajorgrids,
            grid style = dashed,
            width = 8.5cm,
            height = 6cm,
            legend pos = north east
        ]
            \addplot+[very thick, color = blue] table [x = t, y = y, col sep = comma] {data/Cowtemp_Best_reg_coeff_0.01.csv};
            \addplot+[very thick, color = green, mark = o] table [x = t, y = y_rec, col sep = comma] {data/Cowtemp_Best_reg_coeff_0.01.csv};
            \legend{\(y\), \(\hat{y}\)}
        \end{axis}
    \end{tikzpicture}
    \captionof{figure}{The original cowtemp time-series (\(y\)) and the smoothened completion (\(\hat{y}\)) with the forecast.}
    \label{fig:CowTemp}
\end{figure}
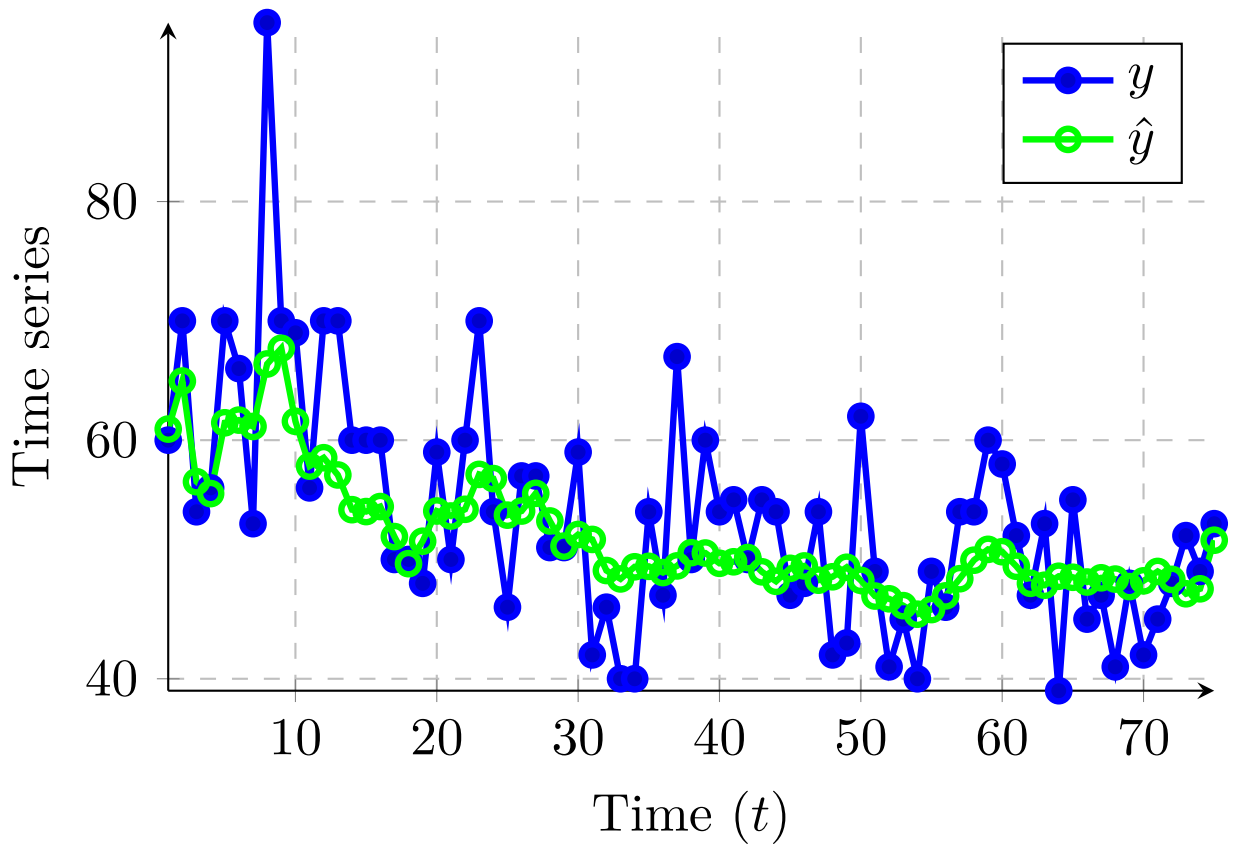

We take the cowtemp time-series, which consists of \(75\) entries and mask the last \(14\) entries with zeros with the rest of the \(61\) entries unchanged. This masked time-series is used to construct a Hankel matrix which is then approximated with rank \(r = 28\) as in \cite{gillard-usevich}. Here, we consider the approximate matrix completion as described in problem~\ref{eqn:MatrixCompletionAlternative2} with the addition of the regularization term \(\frac{\lambda_\mathrm{reg}}{2}\norm[\Fr]{Z - \widetilde{Y}}^2\) in the objective for a regularization coefficient \(\lambda_\mathrm{reg} > 0\), where \(\widetilde{Y}\) is some other matrix that we expect the Hankel matrix generated by the time-series to be close to, which might as well be the zero matrix. We take \(\lambda_\mathrm{reg} = 0.01\). Furthermore, we alternate between this algorithm and a few iterations of a simple alternating projections algorithm \(5\) times, which was found to drastically speed up convergence. We also fix all entries in \(\widetilde{Y}\) to be the median of the last \(\SI{25}{\%}\) of the known entries, which was found to be \(49.0\). The results are shown in \cref{fig:CowTemp}. The RMSE of the forecast was found to be approximately \(4.59596955\).

We consider the daily close prices of the Google stock from August 19, 2004 to November 4, 2026, NVIDIA stock from January 22, 1999 to November 4, 2026, Exxon stock from January 2, 1962 to November 4, 2026, and Coca-Cola stock from January 2, 1962 to November 4, 2026 and draw uniformly randomly \(1000\) windows of length \(8\), and compute the log returns of the stock prices, which is of length \(7\). We mask the last \(\nu\in\{1, 2, 3\}\) entries, form a Hankel matrix and perform approximate matrix completion with \(\lambda_\mathrm{reg} = 10^{-5}\) to forecast the \(\nu\) entries using the known \(7 - \nu\) entries. To determine the rank, we take the minimum rank that captures \(\SI{95}{\%}\) of the squared singular values. We run \(2\) loops of our algorithm with a maximum of \(10\) iterations. A summary of the results are shown in \cref{tab:StocksCorrelation} where \(\mathrm{Pr}\{\rho > 0.5\}\) denotes the proportion of cases for which \(\rho > 0.5\) holds true.

\begin{table}[H]
    \centering
    \captionof{table}{The statistics of the correlation values (\(\rho\)) for different stocks and prediction horizons (\(\nu\)).}
    \begin{tabular}{c c c c c c}
        \toprule
        Stock & \(\nu\) & Mean & Median & SE & \(\mathrm{Pr}\{\rho > 0.5\}\)\\
        \midrule
        Google & \(1\) & \(0.88696\) & \(0.95629\) & \(0.19503\) & \(0.952\)\\
        Google & \(2\) & \(0.80878\) & \(0.88932\) & \(0.22108\) & \(0.905\)\\
        Google & \(3\) & \(0.73063\) & \(0.80663\) & \(0.24075\) & \(0.840\)\\
        NVIDIA & \(1\) & \(0.89537\) & \(0.95202\) & \(0.16795\) & \(0.968\)\\
        NVIDIA & \(2\) & \(0.82056\) & \(0.89512\) & \(0.21653\) & \(0.927\)\\
        NVIDIA & \(3\) & \(0.71947\) & \(0.78212\) & \(0.23852\) & \(0.838\)\\
        Exxon & \(1\) & \(0.87873\) & \(0.95246\) & \(0.20339\) & \(0.949\)\\
        Exxon & \(2\) & \(0.81036\) & \(0.89000\) & \(0.21950\) & \(0.911\)\\
        Exxon & \(3\) & \(0.72771\) & \(0.78818\) & \(0.23248\) & \(0.863\)\\
        Coca-Cola & \(1\) & \(0.88565\) & \(0.95825\) & \(0.20396\) & \(0.953\)\\
        Coca-Cola & \(2\) & \(0.81428\) & \(0.89689\) & \(0.23074\) & \(0.915\)\\
        Coca-Cola & \(3\) & \(0.72664\) & \(0.80079\) & \(0.23729\) & \(0.851\)\\
        \bottomrule
    \end{tabular}
    \label{tab:StocksCorrelation}
\end{table}

\begin{figure}[H]
    \centering
    \begin{tikzpicture}
    \begin{groupplot}[
        group style = {
            group size = 2 by 3,
            horizontal sep = 1.5cm,
            vertical sep = 1.75cm,
        },
        auto title,
    ]
    \nextgroupplot[
        ybar,
        grid = both,
        grid style = {dashed},
        minor grid style = {dashed},
        width = 8cm,
        height = 3cm,
        xlabel = {Rank (\(r\))},
        ylabel = {Frequency},
        xtick = {1, 2, 3, 4},
    ]
        \addplot+[hist = {bins = 4}] table [col sep = comma, y = rank_vals] {./data/Stock_Google_Prediction_Horizon_1.csv};
    \nextgroupplot[
        ybar,
        grid = both,
        grid style = {dashed},
        minor grid style = {dashed},
        width = 8cm,
        height = 3cm,
        xlabel = {Correlation (\(\rho\))},
        ylabel = {Frequency},
    ]
        \addplot+[hist = {bins = 10}] table [col sep = comma, y = rho_vals] {./data/Stock_Google_Prediction_Horizon_1.csv};
    \nextgroupplot[
        ybar,
        grid = both,
        grid style = {dashed},
        minor grid style = {dashed},
        width = 8cm,
        height = 3cm,
        xlabel = {Rank (\(r\))},
        ylabel = {Frequency},
        xtick = {1, 2, 3, 4},
    ]
        \addplot+[hist = {bins = 4}] table [col sep = comma, y = rank_vals] {./data/Stock_Google_Prediction_Horizon_2.csv};
    \nextgroupplot[
        ybar,
        grid = both,
        grid style = {dashed},
        minor grid style = {dashed},
        width = 8cm,
        height = 3cm,
        xlabel = {Correlation (\(\rho\))},
        ylabel = {Frequency},
    ]
        \addplot+[hist = {bins = 10}] table [col sep = comma, y = rho_vals] {./data/Stock_Google_Prediction_Horizon_2.csv};
    \nextgroupplot[
        ybar,
        grid = both,
        grid style = {dashed},
        minor grid style = {dashed},
        width = 8cm,
        height = 3cm,
        xlabel = {Rank (\(r\))},
        ylabel = {Frequency},
        xtick = {1, 2, 3, 4},
    ]
        \addplot+[hist = {bins = 4}] table [col sep = comma, y = rank_vals] {./data/Stock_Google_Prediction_Horizon_3.csv};
    \nextgroupplot[
        ybar,
        grid = both,
        grid style = {dashed},
        minor grid style = {dashed},
        width = 8cm,
        height = 3cm,
        xlabel = {Correlation (\(\rho\))},
        ylabel = {Frequency},
    ]
        \addplot+[hist = {bins = 10}] table [col sep = comma, y = rho_vals] {./data/Stock_Google_Prediction_Horizon_3.csv};
    \end{groupplot}
    \end{tikzpicture}
    \caption{Histogram of the rank and correlation values between the true and approximated time-series for a prediction horizon of \(1\) (a, b), \(2\) (c, d), and \(3\) (e, f) for Google stock.}
    \label{fig:GoogleStockRankCorr}
\end{figure}
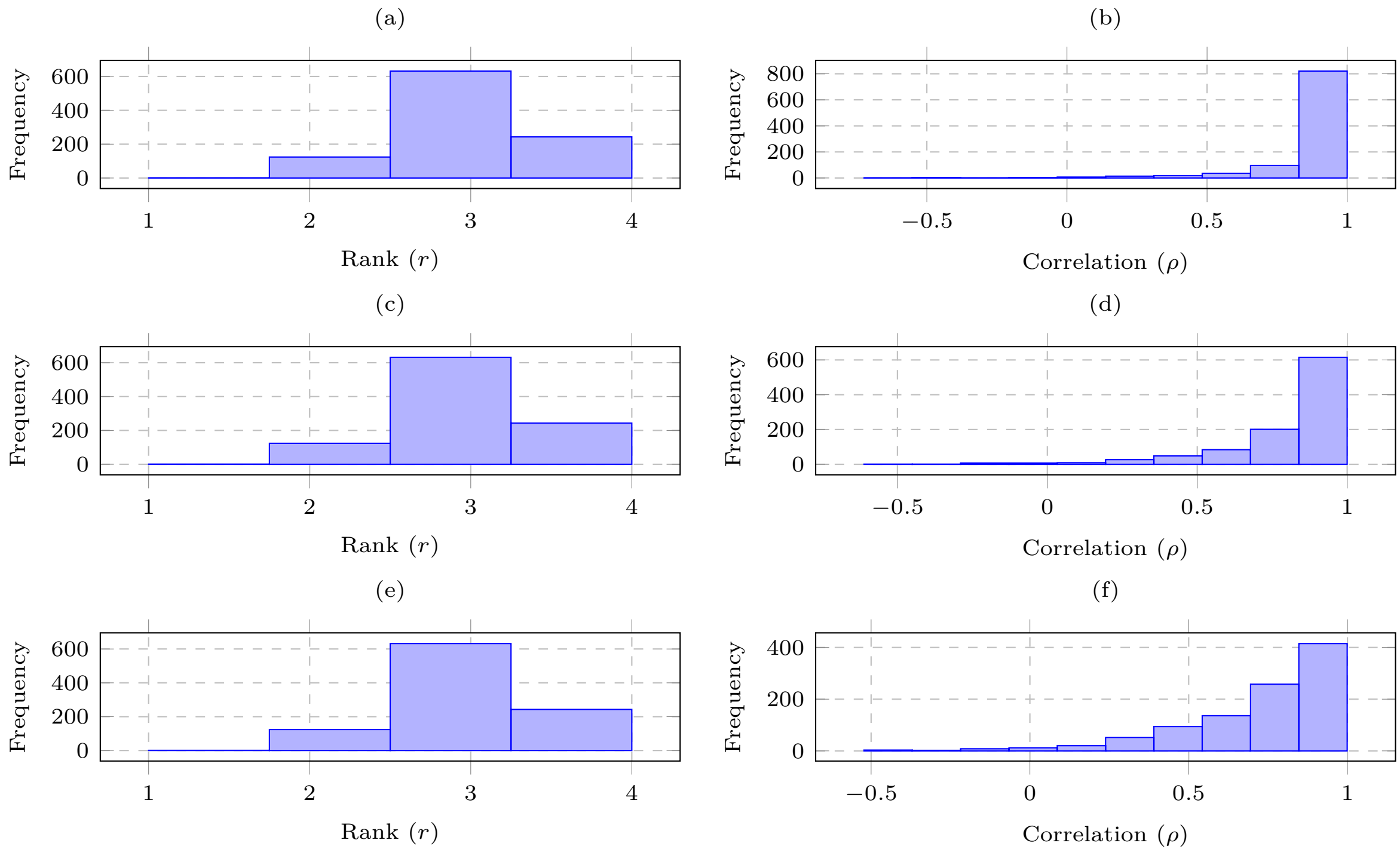

\begin{figure}[H]
    \centering
    \begin{tikzpicture}
    \begin{groupplot}[
        group style = {
            group size = 2 by 3,
            horizontal sep = 1.5cm,
            vertical sep = 1.75cm,
        },
        auto title,
    ]
    \nextgroupplot[
        ybar,
        grid = both,
        grid style = {dashed},
        minor grid style = {dashed},
        width = 8cm,
        height = 3cm,
        xlabel = {Rank (\(r\))},
        ylabel = {Frequency},
        xtick = {1, 2, 3, 4},
    ]
        \addplot+[hist = {bins = 4}] table [col sep = comma, y = rank_vals] {./data/Stock_NVIDIA_Prediction_Horizon_1.csv};
    \nextgroupplot[
        ybar,
        grid = both,
        grid style = {dashed},
        minor grid style = {dashed},
        width = 8cm,
        height = 3cm,
        xlabel = {Correlation (\(\rho\))},
        ylabel = {Frequency},
    ]
        \addplot+[hist = {bins = 10}] table [col sep = comma, y = rho_vals] {./data/Stock_NVIDIA_Prediction_Horizon_1.csv};
    \nextgroupplot[
        ybar,
        grid = both,
        grid style = {dashed},
        minor grid style = {dashed},
        width = 8cm,
        height = 3cm,
        xlabel = {Rank (\(r\))},
        ylabel = {Frequency},
        xtick = {1, 2, 3, 4},
    ]
        \addplot+[hist = {bins = 4}] table [col sep = comma, y = rank_vals] {./data/Stock_NVIDIA_Prediction_Horizon_2.csv};
    \nextgroupplot[
        ybar,
        grid = both,
        grid style = {dashed},
        minor grid style = {dashed},
        width = 8cm,
        height = 3cm,
        xlabel = {Correlation (\(\rho\))},
        ylabel = {Frequency},
    ]
        \addplot+[hist = {bins = 10}] table [col sep = comma, y = rho_vals] {./data/Stock_NVIDIA_Prediction_Horizon_2.csv};
    \nextgroupplot[
        ybar,
        grid = both,
        grid style = {dashed},
        minor grid style = {dashed},
        width = 8cm,
        height = 3cm,
        xlabel = {Rank (\(r\))},
        ylabel = {Frequency},
        xtick = {1, 2, 3, 4},
    ]
        \addplot+[hist = {bins = 4}] table [col sep = comma, y = rank_vals] {./data/Stock_NVIDIA_Prediction_Horizon_3.csv};
    \nextgroupplot[
        ybar,
        grid = both,
        grid style = {dashed},
        minor grid style = {dashed},
        width = 8cm,
        height = 3cm,
        xlabel = {Correlation (\(\rho\))},
        ylabel = {Frequency},
    ]
        \addplot+[hist = {bins = 10}] table [col sep = comma, y = rho_vals] {./data/Stock_NVIDIA_Prediction_Horizon_3.csv};
    \end{groupplot}
    \end{tikzpicture}
    \caption{Histogram of the rank and correlation values between the true and approximated time-series for a prediction horizon of \(1\) (a, b), \(2\) (c, d), and \(3\) (e, f) for NVIDIA stock.}
    \label{fig:NVIDIAStockRankCorr}
\end{figure}
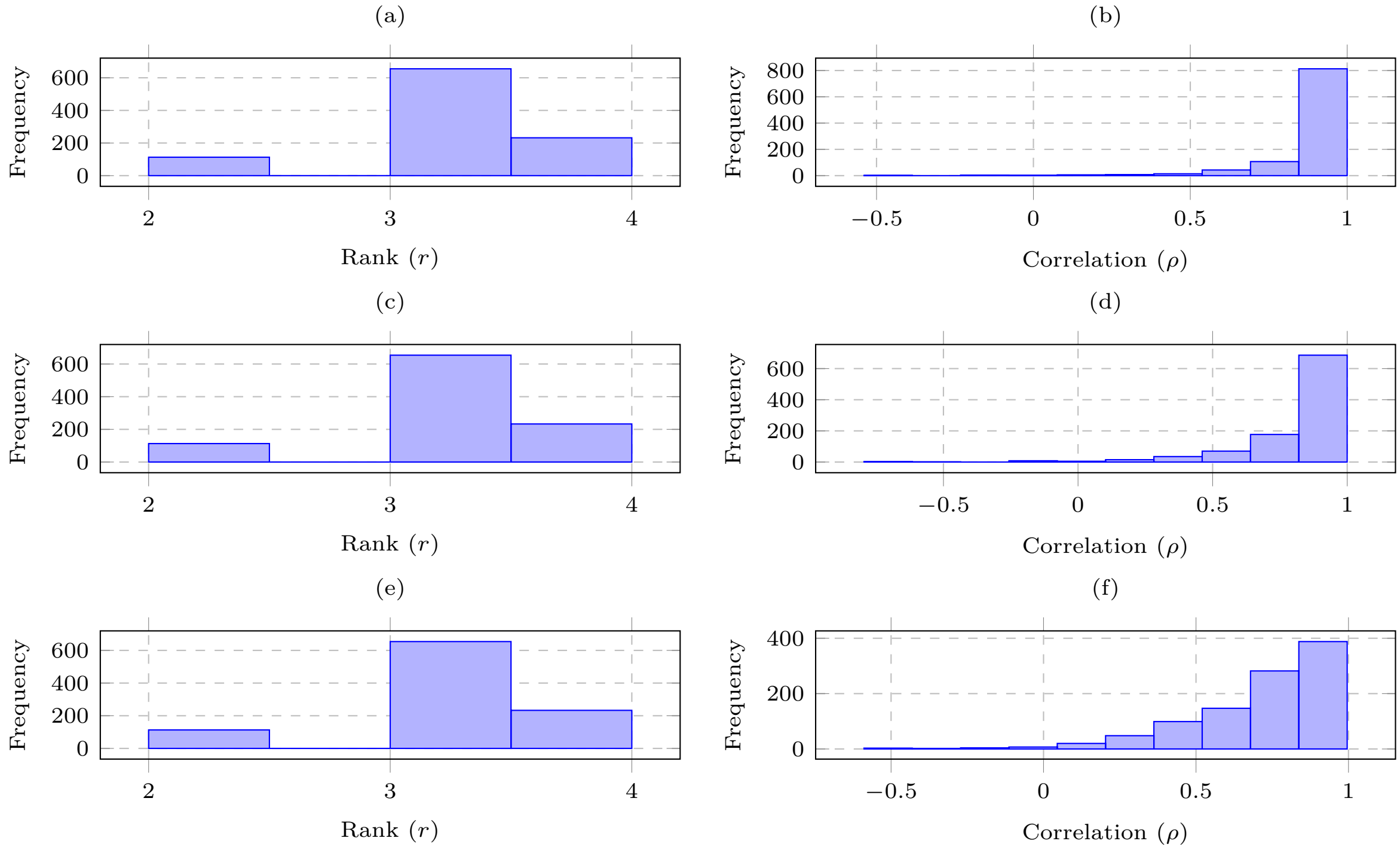

\begin{figure}[H]
    \centering
    \begin{tikzpicture}
    \begin{groupplot}[
        group style = {
            group size = 2 by 3,
            horizontal sep = 1.5cm,
            vertical sep = 1.75cm,
        },
        auto title,
    ]
    \nextgroupplot[
        ybar,
        grid = both,
        grid style = {dashed},
        minor grid style = {dashed},
        width = 8cm,
        height = 3cm,
        xlabel = {Rank (\(r\))},
        ylabel = {Frequency},
        xtick = {1, 2, 3, 4},
    ]
        \addplot+[hist = {bins = 4}] table [col sep = comma, y = rank_vals] {./data/Stock_Exxon_Prediction_Horizon_1.csv};
    \nextgroupplot[
        ybar,
        grid = both,
        grid style = {dashed},
        minor grid style = {dashed},
        width = 8cm,
        height = 3cm,
        xlabel = {Correlation (\(\rho\))},
        ylabel = {Frequency},
    ]
        \addplot+[hist = {bins = 10}] table [col sep = comma, y = rho_vals] {./data/Stock_Exxon_Prediction_Horizon_1.csv};
    \nextgroupplot[
        ybar,
        grid = both,
        grid style = {dashed},
        minor grid style = {dashed},
        width = 8cm,
        height = 3cm,
        xlabel = {Rank (\(r\))},
        ylabel = {Frequency},
        xtick = {1, 2, 3, 4},
    ]
        \addplot+[hist = {bins = 4}] table [col sep = comma, y = rank_vals] {./data/Stock_Exxon_Prediction_Horizon_2.csv};
    \nextgroupplot[
        ybar,
        grid = both,
        grid style = {dashed},
        minor grid style = {dashed},
        width = 8cm,
        height = 3cm,
        xlabel = {Correlation (\(\rho\))},
        ylabel = {Frequency},
    ]
        \addplot+[hist = {bins = 10}] table [col sep = comma, y = rho_vals] {./data/Stock_Exxon_Prediction_Horizon_2.csv};
    \nextgroupplot[
        ybar,
        grid = both,
        grid style = {dashed},
        minor grid style = {dashed},
        width = 8cm,
        height = 3cm,
        xlabel = {Rank (\(r\))},
        ylabel = {Frequency},
        xtick = {1, 2, 3, 4},
    ]
        \addplot+[hist = {bins = 4}] table [col sep = comma, y = rank_vals] {./data/Stock_Exxon_Prediction_Horizon_3.csv};
    \nextgroupplot[
        ybar,
        grid = both,
        grid style = {dashed},
        minor grid style = {dashed},
        width = 8cm,
        height = 3cm,
        xlabel = {Correlation (\(\rho\))},
        ylabel = {Frequency},
    ]
        \addplot+[hist = {bins = 10}] table [col sep = comma, y = rho_vals] {./data/Stock_Exxon_Prediction_Horizon_3.csv};
    \end{groupplot}
    \end{tikzpicture}
    \caption{Histogram of the rank and correlation values between the true and approximated time-series for a prediction horizon of \(1\) (a, b), \(2\) (c, d), and \(3\) (e, f) for Exxon stock.}
    \label{fig:ExxonStockRankCorr}
\end{figure}
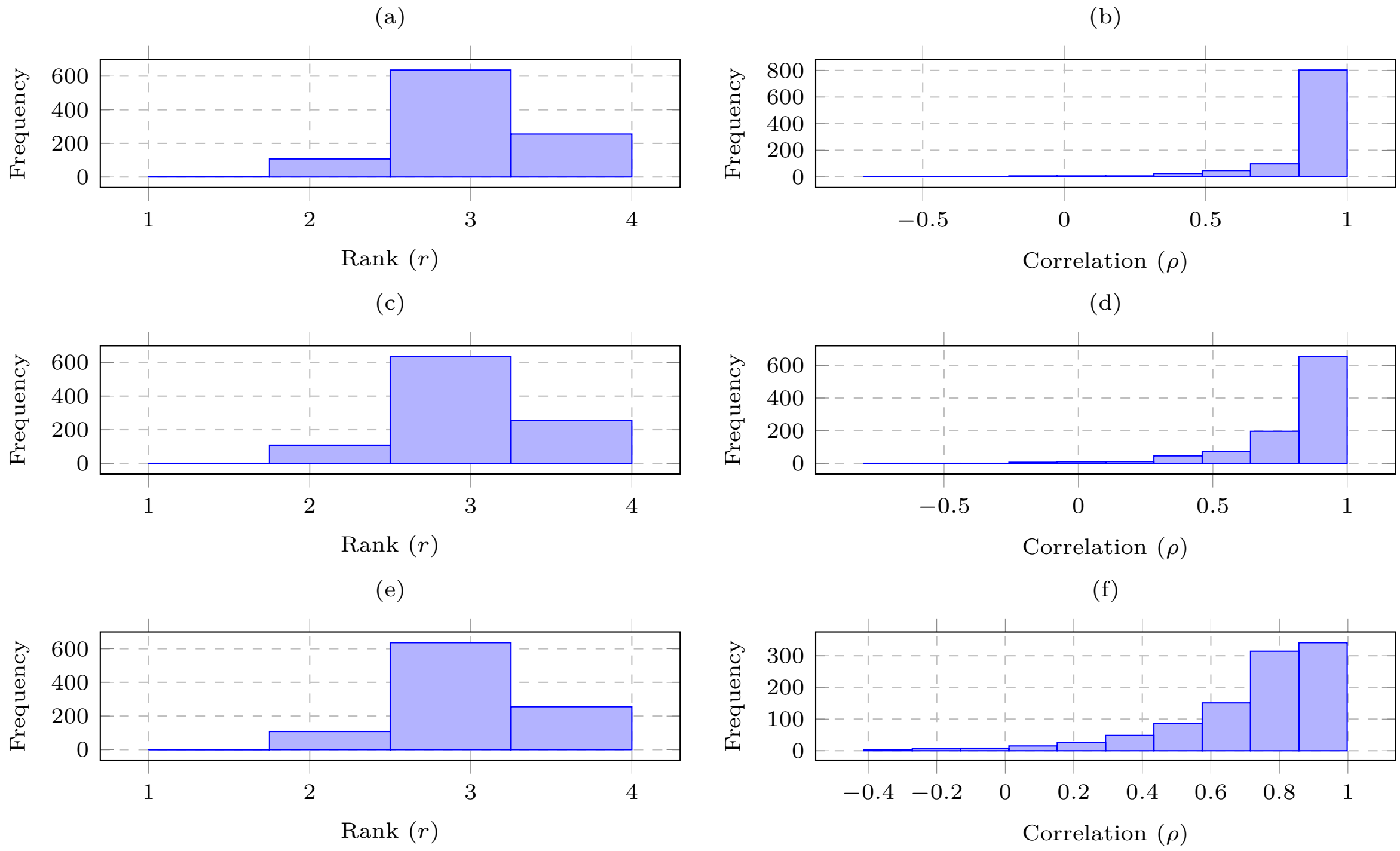

\begin{figure}[H]
    \centering
    \begin{tikzpicture}
    \begin{groupplot}[
        group style = {
            group size = 2 by 3,
            horizontal sep = 1.5cm,
            vertical sep = 1.75cm,
        },
        auto title,
    ]
    \nextgroupplot[
        ybar,
        grid = both,
        grid style = {dashed},
        minor grid style = {dashed},
        width = 8cm,
        height = 3cm,
        xlabel = {Rank (\(r\))},
        ylabel = {Frequency},
        xtick = {1, 2, 3, 4},
    ]
        \addplot+[hist = {bins = 4}] table [col sep = comma, y = rank_vals] {./data/Stock_Coca_Cola_Prediction_Horizon_1.csv};
    \nextgroupplot[
        ybar,
        grid = both,
        grid style = {dashed},
        minor grid style = {dashed},
        width = 8cm,
        height = 3cm,
        xlabel = {Correlation (\(\rho\))},
        ylabel = {Frequency},
    ]
        \addplot+[hist = {bins = 10}] table [col sep = comma, y = rho_vals] {./data/Stock_Coca_Cola_Prediction_Horizon_1.csv};
    \nextgroupplot[
        ybar,
        grid = both,
        grid style = {dashed},
        minor grid style = {dashed},
        width = 8cm,
        height = 3cm,
        xlabel = {Rank (\(r\))},
        ylabel = {Frequency},
        xtick = {1, 2, 3, 4},
    ]
        \addplot+[hist = {bins = 4}] table [col sep = comma, y = rank_vals] {./data/Stock_Coca_Cola_Prediction_Horizon_2.csv};
    \nextgroupplot[
        ybar,
        grid = both,
        grid style = {dashed},
        minor grid style = {dashed},
        width = 8cm,
        height = 3cm,
        xlabel = {Correlation (\(\rho\))},
        ylabel = {Frequency},
    ]
        \addplot+[hist = {bins = 10}] table [col sep = comma, y = rho_vals] {./data/Stock_Coca_Cola_Prediction_Horizon_2.csv};
    \nextgroupplot[
        ybar,
        grid = both,
        grid style = {dashed},
        minor grid style = {dashed},
        width = 8cm,
        height = 3cm,
        xlabel = {Rank (\(r\))},
        ylabel = {Frequency},
        xtick = {1, 2, 3, 4},
    ]
        \addplot+[hist = {bins = 4}] table [col sep = comma, y = rank_vals] {./data/Stock_Coca_Cola_Prediction_Horizon_3.csv};
    \nextgroupplot[
        ybar,
        grid = both,
        grid style = {dashed},
        minor grid style = {dashed},
        width = 8cm,
        height = 3cm,
        xlabel = {Correlation (\(\rho\))},
        ylabel = {Frequency},
    ]
        \addplot+[hist = {bins = 10}] table [col sep = comma, y = rho_vals] {./data/Stock_Coca_Cola_Prediction_Horizon_3.csv};
    \end{groupplot}
    \end{tikzpicture}
    \caption{Histogram of the rank and correlation values between the true and approximated time-series for a prediction horizon of \(1\) (a, b), \(2\) (c, d), and \(3\) (e, f) for Coca-Cola stock.}
    \label{fig:CocaColaStockRankCorr}
\end{figure}
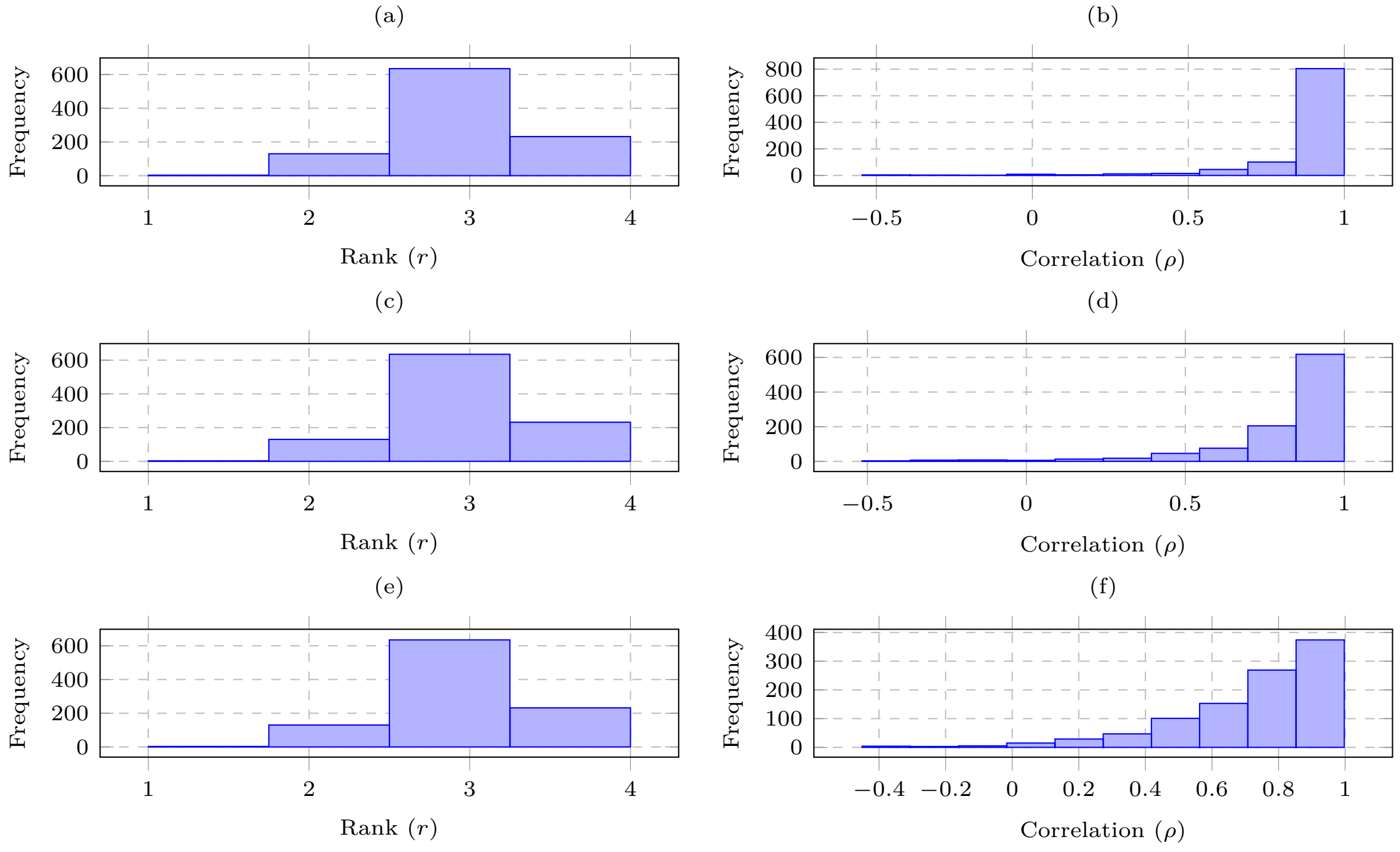

\subsubsection{Asymptotic Complexity}
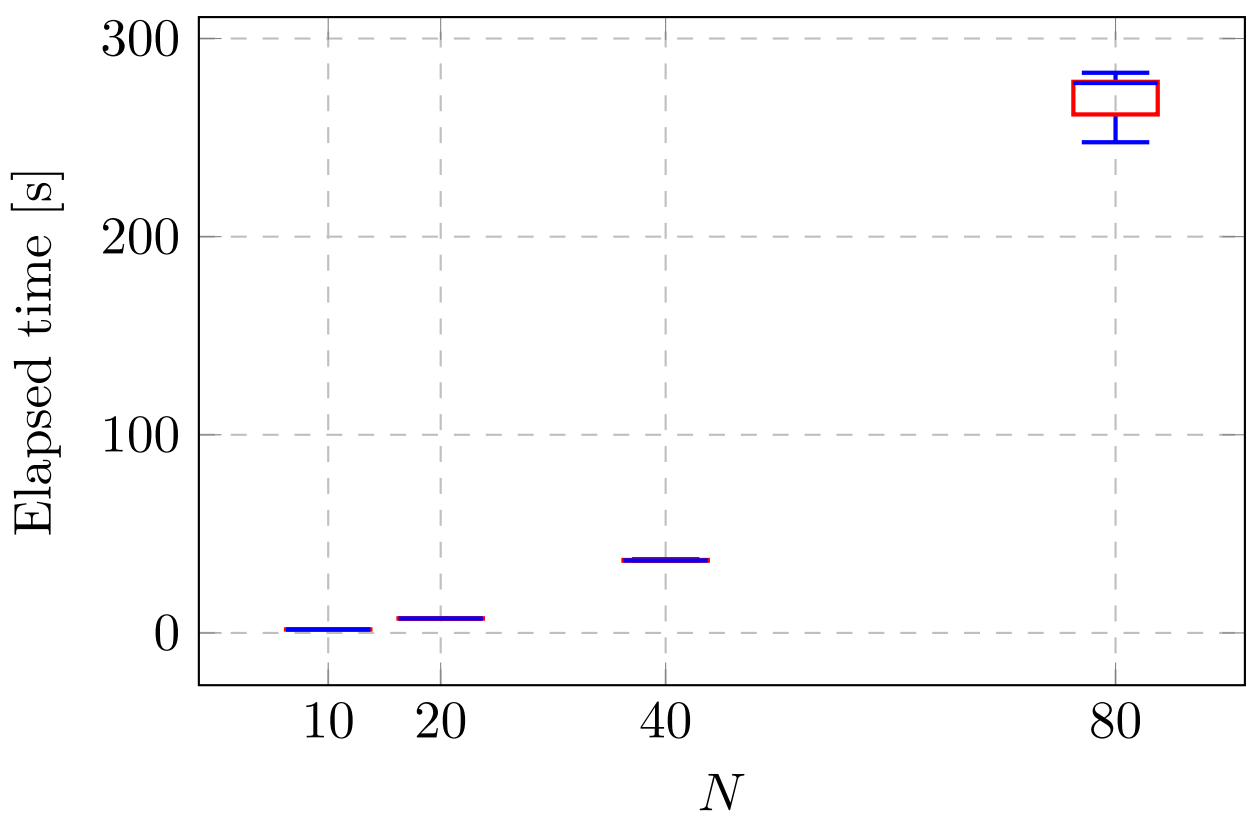
\begin{figure}[H]
    \centering
    \begin{tikzpicture}
    \begin{axis}[
        boxplot/every box/.style = {thick, solid, red},
        boxplot/every median/.style = {thick, solid, blue},
        boxplot/every whisker/.style = {thick, solid, blue},
        boxplot/draw direction = y,
        grid = both,
        grid style = {dashed},
        minor grid style = {dashed},
        width = 8.5cm,
        height = 6cm,
        xlabel = {\(N\)},
        ylabel = {Elapsed time [\(\si{\second}\)]},
        xtick = {10, 20, 40, 80}
    ]
        \addplot+ [boxplot prepared = {
            draw position = 10.00000,
            lower whisker = 1.72092, lower quartile = 1.72092,
            median = 1.75749,
            upper quartile = 1.78048, upper whisker = 1.91621,
            box extend = 7.50000
            },
        ] coordinates {};
        \addplot+ [boxplot prepared = {
            draw position = 20.00000,
            lower whisker = 7.21687, lower quartile = 7.24797,
            median = 7.30166,
            upper quartile = 7.37870, upper whisker = 7.43415,
            box extend = 7.50000
            },
        ] coordinates {};
        \addplot+ [boxplot prepared = {
            draw position = 40.00000,
            lower whisker = 36.39163, lower quartile = 36.46907,
            median = 36.59862,
            upper quartile = 36.82090, upper whisker = 37.19165,
            box extend = 7.50000
            },
        ] coordinates {};
        \addplot+ [boxplot prepared = {
            draw position = 80.00000,
            lower whisker = 247.64866, lower quartile = 261.68583,
            median = 277.54471,
            upper quartile = 278.14804, upper whisker = 282.71779,
            box extend = 7.50000
            },
        ] coordinates {};
    \end{axis}
    \end{tikzpicture}
    \caption{Time elapsed for \(10\) iterations of the algorithm for completing a \(N\times N\) Hankel matrix with \(p = 0.5\).}
    \label{fig:StructuredCompletionComplexity}
\end{figure}

We take the weather dataset \cite{weather} and construct a Hankel matrix by progressively selecting more data, mask the entries with probability of omission \(p = 0.5\), and try to complete it with rank \(r = 5\). The algorithm is run for \(10\) iterations \(5\) times for each \(N\in\{10, 20, 40, 80\}\), where the Hankel matrix formed is of size \(N\times N\). The results are given in \cref{fig:StructuredCompletionComplexity}. From a least-squares fit after taking the logarithm of the elapsed times, the elapsed time seems to grow as \(\Order\left(N^{2.5}\right)\).

\section{Conclusion}
We have derived the constrained critical points of a generalized version of a low-rank approximation with affine equality constraints problem and have proposed an algorithm to find a ``good'' local minimum of the problem. Furthermore, we have demonstrated the proposed algorithm on a variety of applications with comparisons against some standard benchmark algorithms. We also note that inequality constraints can easily be incorporated into the problem through penalty functions, especially when they are given by affine functions as the Hessian in \cref{thm:CriticalPoints} remains unchanged.
\printbibliography

@article{haeffele-vidal-2020,
  author={Benjamin D. Haeffele and  René Vidal},
  journal={IEEE Transactions on Pattern Analysis and Machine Intelligence}, 
  title={Structured Low-Rank Matrix Factorization: Global Optimality, Algorithms, and Applications}, 
  year={2020},
  volume={42},
  number={6},
  pages={1468-1482}
}

@book{bertsekas,
  author = {Dimitri P. Bertsekas},
  title = {Nonlinear Programming},
  publisher = {Athena Scientific},
  edition = {Second Edition},
  year = 2016
}

@misc{ishteva-usevich-markovsky,
  title={Factorization approach to structured low-rank approximation with applications}, 
  author={Mariya Ishteva and Konstantin Usevich and Ivan Markovsky},
  year={2014},
  eprint={1308.1827},
  archivePrefix={arXiv},
  primaryClass={math.NA},
  url={https://arxiv.org/abs/1308.1827}, 
}

@misc{gillard-usevich,
  title={Hankel low-rank approximation and completion in time series analysis and forecasting: a brief review}, 
  author={Jonathan Gillard and Konstantin Usevich},
  year={2022},
  eprint={2206.05103},
  archivePrefix={arXiv},
  primaryClass={math.NA},
  url={https://arxiv.org/abs/2206.05103}, 
}

@article{cifuentes,
  author = {Cifuentes, Diego},
  title = {A Convex Relaxation to Compute the Nearest Structured Rank Deficient Matrix},
  journal = {SIAM Journal on Matrix Analysis and Applications},
  volume = {42},
  number = {2},
  pages = {708-729},
  year = {2021}
}

@article{gillard-zhigljavsky,
  author = {Jonathan W. Gillard and Anatoly A. Zhigljavsky},
  title = {Stochastic algorithms for solving structured low-rank matrix approximation problems},
  journal = {Communications in Nonlinear Science and Numerical Simulation},
  volume = {21},
  number = {1},
  pages = {70-88},
  year = {2015},
  note = {Numerical Computations: Theory and Algorithms (NUMTA 2013), International Conference and Summer School},
  issn = {1007-5704}
}

@article{zhi-approximate-gcd,
  title={Computing approximate GCD of univariate polynomials by structure total least norm},
  author={Zhi, Lihong and Yang, Zhengfeng},
  journal={MM research preprints},
  volume={24},
  pages={375--387},
  year={2004}
}

@article{cadzow,
  author={Cadzow, J.A.},
  journal={IEEE Transactions on Acoustics, Speech, and Signal Processing}, 
  title={Signal enhancement-a composite property mapping algorithm}, 
  year={1988},
  volume={36},
  number={1},
  pages={49-62},
  doi={10.1109/29.1488}
}

@article{knirsch-petz-plonka,
  title = {Optimal rank-1 Hankel approximation of matrices: Frobenius norm and spectral norm and Cadzow's algorithm},
  journal = {Linear Algebra and its Applications},
  volume = {629},
  pages = {1-39},
  year = {2021},
  issn = {0024-3795},
  doi = {https://doi.org/10.1016/j.laa.2021.07.004},
  url = {https://www.sciencedirect.com/science/article/pii/S0024379521002627},
  author = {Hanna Knirsch and Markus Petz and Gerlind Plonka}
}

@article{gillard-zhigljavsky-2013,
  author = {Gillard, Jonathan and Zhigljavsky, Anatoly},
  title = {Optimization challenges in the structured low rank approximation problem},
  journal = {Journal of Global Optimization},
  year = {2013},
  pages = {733-751},
  volume = {57},
  issue = {3}
}

@misc{zvonarev-golyandina,
  title={Iterative algorithms for weighted and unweighted finite-rank time-series approximations}, 
  author={Nikita Zvonarev and Nina Golyandina},
  year={2015},
  eprint={1507.02751},
  archivePrefix={arXiv},
  primaryClass={stat.ME},
  url={https://arxiv.org/abs/1507.02751}, 
}

@article{markovsky-usevich,
  title = {Software for weighted structured low-rank approximation},
  journal = {Journal of Computational and Applied Mathematics},
  volume = {256},
  pages = {278-292},
  year = {2014},
  issn = {0377-0427},
  doi = {https://doi.org/10.1016/j.cam.2013.07.048},
  author = {Ivan Markovsky and Konstantin Usevich}
}

@article{recht-fazel-parrilo,
  author = {Recht, Benjamin and Fazel, Maryam and Parrilo, Pablo A.},
  title = {Guaranteed Minimum-Rank Solutions of Linear Matrix Equations via Nuclear Norm Minimization},
  journal = {SIAM Review},
  volume = {52},
  number = {3},
  pages = {471-501},
  year = {2010},
  doi = {10.1137/070697835}
}

@online{traffic,
  title = {Traffic Dataset},
  url = {https://github.com/laiguokun/multivariate-time-series-data},
}

@online{weather,
  title = {Weather Dataset},
  note = {
  \url{https://www.bgc-jena.mpg.de/wetter/},
  \url{https://www.kaggle.com/datasets/giochelavaipiatti/time-series-forecasts-popular-benchmark-datasets?select=weather.csv}
  },
}

@article{coulson-2021,
  title={Distributionally robust chance constrained data-enabled predictive control},
  author={Coulson, Jeremy and Lygeros, John and D{\"o}rfler, Florian},
  journal={IEEE Transactions on Automatic Control},
  volume={67},
  number={7},
  pages={3289--3304},
  year={2021},
  publisher={IEEE}
}

@article{esfahani-kuhn,
  author = {Mohajerin Esfahani, Peyman and Kuhn, Daniel},
  journal = {Mathematical Programming},
  volume = {171},
  issue = {1},
  title = {Data-driven distributionally robust optimization using the Wasserstein metric: performance guarantees and tractable reformulations},
  year = {2018},
  pages = {115-166},
}

@article{candes-plan-yaniv,
  author={Candès, Emmanuel J. and Plan, Yaniv},
  journal={IEEE Transactions on Information Theory}, 
  title={Tight Oracle Inequalities for Low-Rank Matrix Recovery From a Minimal Number of Noisy Random Measurements}, 
  year={2011},
  volume={57},
  number={4},
  pages={2342-2359}
}

@article{candes-tao,
  author={Candes, Emmanuel J. and Tao, Terence},
  journal={IEEE Transactions on Information Theory}, 
  title={The Power of Convex Relaxation: Near-Optimal Matrix Completion}, 
  year={2010},
  volume={56},
  number={5},
  pages={2053-2080}
}

@article{candes-recht,
  author = {Cand{\`e}s, Emmanuel and Recht, Benjamin},
  title = {Exact matrix completion via convex optimization},
  year = {2012},
  issue_date = {June 2012},
  publisher = {Association for Computing Machinery},
  address = {New York, NY, USA},
  volume = {55},
  number = {6},
  issn = {0001-0782},
  journal = {Commun. ACM},
  month = jun,
  pages = {111-119},
  numpages = {9}
}

@article{vandereycken,
  author = {Vandereycken, Bart},
  title = {Low-Rank Matrix Completion by Riemannian Optimization},
  journal = {SIAM Journal on Optimization},
  volume = {23},
  number = {2},
  pages = {1214-1236},
  year = {2013}
}

@ARTICLE{candes-plan,
  author={Candes, Emmanuel J. and Plan, Yaniv},
  journal={Proceedings of the IEEE}, 
  title={Matrix Completion With Noise}, 
  year={2010},
  volume={98},
  number={6},
  pages={925-936}
}

@article{negahban-wainwright-2012,
  title={Restricted strong convexity and weighted matrix completion: Optimal bounds with noise},
  author={Negahban, Sahand and Wainwright, Martin J},
  journal={The Journal of Machine Learning Research},
  volume={13},
  number={1},
  pages={1665--1697},
  year={2012},
  publisher={JMLR. org}
}

@article{chen-chi-fan-ma-yan,
  title={Noisy matrix completion: Understanding statistical guarantees for convex relaxation via nonconvex optimization},
  author={Chen, Yuxin and Chi, Yuejie and Fan, Jianqing and Ma, Cong and Yan, Yuling},
  journal={SIAM journal on optimization},
  volume={30},
  number={4},
  pages={3098--3121},
  year={2020},
  publisher={SIAM}
}

@article{negahban-wainwright-2011,
  author = {Sahand Negahban and Martin J. Wainwright},
  title = {{Estimation of (near) low-rank matrices with noise and high-dimensional scaling}},
  volume = {39},
  journal = {The Annals of Statistics},
  number = {2},
  publisher = {Institute of Mathematical Statistics},
  pages = {1069 -- 1097},
  year = {2011}
}

@article{fazel-pong-sun-tseng,
  author = {Fazel, Maryam and Pong, Ting Kei and Sun, Defeng and Tseng, Paul},
  title = {Hankel Matrix Rank Minimization with Applications to System Identification and Realization},
  journal = {SIAM Journal on Matrix Analysis and Applications},
  volume = {34},
  number = {3},
  pages = {946-977},
  year = {2013}
}

@article{zhang-zhang,
  author={Zhang, Shuimei and Zhang, Yimin D.},
  journal={IEEE Transactions on Signal Processing}, 
  title={Low-Rank Hankel Matrix Completion for Robust Time-Frequency Analysis}, 
  year={2020},
  volume={68},
  number={},
  pages={6171-6186}
}

@article{reddi,
  title = {Eigenvector properties of Toeplitz matrices and their application to spectral analysis of time series},
  journal = {Signal Processing},
  volume = {7},
  number = {1},
  pages = {45-56},
  year = {1984},
  issn = {0165-1684},
  author = {S.S Reddi}
}

@ARTICLE{wahlberg,
  author={Wahlberg, B.},
  journal={IEEE Transactions on Automatic Control}, 
  title={System identification using Kautz models}, 
  year={1994},
  volume={39},
  number={6},
  pages={1276-1282},
}

@article{strang,
  author = {Strang, Gilbert},
  title = {A Proposal for Toeplitz Matrix Calculations},
  journal = {Studies in Applied Mathematics},
  volume = {74},
  number = {2},
  pages = {171-176},
  year = {1986}
}

@article{bitmead-kung-anderson-kailath,
  author={Bitmead, R. and Kung, S.-Y. and Anderson, B. and Kailath, T.},
  journal={IEEE Transactions on Automatic Control}, 
  title={Greatest common divisor via generalized Sylvester and Bezout matrices}, 
  year={1978},
  volume={23},
  number={6},
  pages={1043-1047}
}
\end{document}